\documentclass[10pt, reqno]{article}
\usepackage{amsmath,amsthm,amsfonts,amssymb,epsfig}

\newtheorem{tm}{Theorem}[section]
\newtheorem{lemma}[tm]{Lemma}
\newtheorem{prop}[tm]{Proposition}
\newtheorem{cor}[tm]{Corollary}

 \theoremstyle{definition}
 \newtheorem{definition}{Definition}

\newtheorem{alg}{Algorithm}

\newcommand{\beqa}{\begin{eqnarray*}}
\newcommand{\eeqa}{\end{eqnarray*}}

\DeclareMathOperator*{\supp}{supp}

\newcommand{\bH}{\field{H}}        %

 \def\cF{\mathcal{F}}              

 \def\cH{\mathcal{H}}
 \def\cB{\mathcal{B}}

 \def\cA{\mathcal{A}}

 \def\cC{\mathcal{C}}

 \def\cN{\mathcal{N}}

\def\bV{\mathbf{V}}
\def\bl{\mathbf{l}}
\def\ba{\mathbf{A}}
\def\bA1{\mathbf{A}_1}
\def\bH{\mathbf{H}}

\newcommand{\vv}{{\text{\bf v}}}
\newcommand{\uu}{{\text{\bf u}}}
\newcommand{\nn}{{\text{\bf n}}}
\newcommand{\RR}{{\mathbb{R}}}
\newcommand{\LL}{{\boldsymbol{L}}}
\newcommand{\JJ}{{\boldsymbol{J}}}

\newcommand{\BB}{{\boldsymbol{B}}}
\newcommand{\EE}{{\boldsymbol{E}}}
\newcommand{\KK}{{\boldsymbol{K}}}
\newcommand{\g}{{\text{\bf g}}}

 \newcommand{\XX}{{\text{\bf X}}}

\def\L{\left(}
\def\R{\right)}

\def\<{\left<}
\def\>{\right>}

\def\mv1{M_v^1}

\newcommand{\norm}[1]{\lVert#1\rVert}

\newcommand{\NN}{\mathbb N}

\newcommand{\PP}{\mathbf P}

\DeclareMathOperator{\Ran}{ran}
\DeclareMathOperator{\Ker}{ker}
\DeclareMathOperator{\id}{id}

\begin{document}
\title{Adaptive Frame Methods for Magnetohydrodynamic Flows}
\author{M.~Charina, C.~Conti, and M.~Fornasier}
\maketitle

\begin{abstract} In this paper we develop   
adaptive numerical schemes for certain nonlinear variational problems.
The discretization of the variational problems is done by
a suitable frame decomposition of the solution, i.e.,
a complete, stable, and redundant expansion. The discretization
yields an equivalent nonlinear problem   on
$\ell_2(\mathcal{N})$, the space of   frame coefficients. The
discrete problem is then adaptively solved using approximated
nested fixed point and Richardson type iterations. We investigate the 
convergence, 
stability, and optimal complexity of the scheme. 
This constitutes a theoretical advantage, for example, with
respect to adaptive finite element schemes for which convergence
and complexity results are usually hard to prove.
The use of frames is further motivated by their redundancy, which, at 
least numerically, 
has been shown to improve the conditioning of the corresponding discretization
matrices.
Also frames are usually easier to construct than Riesz bases. 
Finally, we show how to apply the adaptive scheme we
propose  for finding an approximation to the solution 
of the PDE governing magnetohydrodynamic (MHD) flows, once 
suitable frames are constructed.
\end{abstract}


\section{Introduction}

Adaptive numerical methods have yielded very promising results
\cite{BR2,BW,Ra1,BEK,CDD1,CDD2,DDHS,Do,Ve2} when applied to a large
class of operator equations, in particular, PDE and integral equations.
In  classical schemes the adaptivity is realized at the level
of the discretization and in the finite element space. The finite element space is refined and
enriched  locally at each iteration step depending on some (a posteriori) error
estimators \cite{Doe,Noch}. A novel paradigm for adaptive schemes has been recently proposed
by Cohen, Dahmen, and DeVore in \cite{CDD1,CDD2}, where the discretization via wavelet decompositions is
fixed at the beginning. The adaptivity is indeed realized at the level
 of the solver of the equivalent bi-infinite system of linear equations.
The basic idea is to transform the original problem of PDE into a discrete
 (bi-infinite) linear problem on  $\ell_2({\mathcal N})$, the space of wavelet
 coefficients. The discrete problem is then solved with the
help of approximate iterative schemes.
The advantage of the latter approach is the fact that its convergence and stability can be proved
and its optimal complexity can be estimated asymptotically in terms of the number 
of algebraic operations needed.
On the contrary, it has been a very hard technical
problem to obtain such nice theoretical
estimates for classical finite element methods, although some important theoretical 
results has recently appeared in \cite{BDD,S2} for linear elliptic equations.
A version of the paradigm in \cite{CDD1,CDD2} has been recently proposed
 also for nonlinear  problems  in \cite{CDD3}. It is again based on
(wavelet) bases discretizations. 

One drawback of the wavelet
approach is the construction of the wavelet system itself
especially on domains with complicated geometry or manifolds \cite{DKU,DS1,DS2}. 
The wavelet bases constructed so far exhibit relatively high
condition numbers and limited smoothness. In particular, the patching used to construct global
smooth wavelets by domain decomposition techniques appears
complicated and, in most cases, makes the conditioning even
worse. In fact the global smoothness of the basis, when
implementing  adaptive schemes in \cite{CDD1,CDD2}, is a
necessary condition for getting compressibility (i.e., finitely
banded approximations) of (bi--infinite) stiffness matrices, especially for high order operators. This
bottleneck has led to generalizations of the Cohen, Dahmen and
DeVore approach. These generalizations are based on {\it frame
discretizations }, i.e., stable, redundant, non-orthogonal
expansions \cite{DFR,S}, which are much more flexible and simpler
to construct even on domains of complicated geometry. Frame
construction is usually implemented by Overlapping Domain
Decompositions (ODD) so that patching at the interfaces is no 
 more needed to obtain global smoothness. 
Moreover, the use of frames, due to their intrinsic redundancy, 
improves the conditioning of the corresponding discretization
matrices.
Certainly, ODD generates regions of the domain  where the side
effect  of the redundancy is that functions are no longer uniquely
representable by the global frame system. At first sight, it may seem that  
redundancy contradicts the minimality requirement on the amount of
information being used to approximate the solution. Especially in
fluid mechanics, accurate simulations already require 
processing of huge amounts of data. How can one attempt such
computations if the information is also made redundant? A figurative
answer to this question is the so-called ``dictionary example'':
The larger and richer is my dictionary the {\it shorter} are
the phrases I compose. The use of proper
terminology avoids long circumlocutions for describing an
object. Of course, the key point is the capability to choose the
right terminology. Back to mathematical terms, the combination of
adaptivity (i.e., the capability to choose the right terminology)
and redundancy (i.e., the richness or non-uniqueness of
representations) indeed gives rise to fast and accurate
approximations  \cite{DMA,GV,MZ,Tem,Tr}. Numerical experiments in \cite{W} show that frames
improve conditioning without increasing the effective dimension of
the problem, i.e. without increasing the number of the  
relevant quantities needed for computations.

This encourages us to present a generalization of the approach proposed in \cite{CDD3} to (wavelet) frame discretizations
for some specific nonlinear PDE, in particular those describing certain magnetohydrodynamics problems.
Magnetohydrodynamics   \cite{CMS,HY,MS1,MS3,MS2,R,SH,WW} studies macroscopic
interactions between magnetic field  and fluid conductors of
electricity. In particular, the following physical phenomena are studied:
A flow of an electrically conducting fluid
across magnetic lines causes an electric current in the fluid. The
electric current alters the electromagnetic state of the system
modifying the total magnetic field, which creates the current.
The flow of electric current across magnetic lines is
associated with a body force - Lorentz force - which influences
the fluid flow. To model the behavior of electrically conducting
fluid, 
the stationary, incompressible Navier-Stokes and Maxwell
equations, coupled via Ohm's law and Lorentz force are being used.
We refer to \cite{MS3,MS2} for a rigorous analysis of a non-adaptive finite--element scheme for the 
simulation of MHD flows confined
to a cubic domain only and arising during electromagnetic purification of molten metals before 
the casting stage.
\\
We present a way for transforming the nonlinear magnetohydrodynamics problem, possibly defined even on more general domains,
into an equivalent nonlinear discrete problem  on $\ell_2({\mathcal N})$ by using suitable frame expansions. 
We  show how the discrete
problem  can be  solved {\it adaptively} by means of nested  fixed point and
approximated Richardson iterations. We also discuss convergence, stability,
and, under certain additional assumptions, computational cost (quasi-optimal complexity) of the proposed adaptive procedure.

The paper is organized as follows:
in Section 2 we recall the mathematical model and the equations 
governing MHD flows, specify the boundary conditions, the solution and 
test function spaces. In Section 3, the corresponding weak 
formulation of the physical problem is derived. Its equivalence  to a variational nonlinear 
problem on a suitable subspace of the solution function space 
follows from  the standard LBB 
(Ladyzhenskaya--Babu\v{s}ka--Brezzi) theory. The proofs of the   
results listed in this section can be found in \cite{CMS,MS1,MS2}. 
In Section 4,  the nonlinear variational problem is re-formulated as an 
equivalent fixed point iteration scheme, 
where, at each iteration step, a linear (non--symmetric) elliptic operator equation 
is to be solved. Next, we present the way for discretizing the 
fixed point iteration associated to an abstract nonlinear problem arising 
in MHD. We translate the original nonlinear variational problem into a problem 
on $\ell_2(\mathcal{N})$, the space of suitable {\it frame} coefficients. The 
concept of frames, i.e., stable, redundant, and complete 
expansions, is recalled in Subsection 4.1. In Subsection 4.2 we show how a 
linear (non-symmetric) elliptic operator equation is discretized by 
means of frame expansions and how Algorithm \ref{alg1} is used to approximate its 
solution adaptively, up to any prescribed accuracy.
Using Algorithm \ref{alg1}  as the main building block, we formulate in Section 5 a fully 
discrete and finite version of the fixed point iteration. We show that this
discrete version of the fixed point iteration 
converges to some frame coefficients of the true solution of the original problem. 
In Section 6 we discuss under which conditions on the building block procedures in Algorithm 
\ref{alg1} the suggested scheme performs quasi-optimally with 
respect to suitable {\it sparseness classes} of frame coefficients.
In particular, we show how the flexibility and redundancy of frames lead to technical 
difficulties, which do not arise in case of Riesz bases, when showing complexity estimates. 
In Section 7, we show that the relevant solution space for the MHD problem is  
a product of suitable divergence--free vector spaces.
We shortly discuss the existence and the construction of suitable divergence-free 
frames for such spaces. This allows for applications of the general adaptive scheme
we propose.

Throughout this paper `$a \sim b$' means that both quantities are uniformly bounded by
some constant multiple of each other. Likewise, `$a \lesssim b$' means that there
exists a positive constant $C$ such that $a \le C b$. We determine the constants
explicitly only if their value is crucial for further analysis.
The symbol $\|\cdot\|$, when applied to bounded operators, denotes the operator norm from 
its domain space to its image space, these are not always explicitly specified 
for notational simplicity.

\section{Equations and Boundary Conditions}

In this section, we   model the stationary flow of a viscous, incompressible, electrically 
conducting fluid  occupying
a 3--dimensional bounded region $\Omega \subset \RR^3$. For later analysis, we assume that
$\Omega$ is a Lipschitz domain. Assume also that ${\bf F}$, a body force, and $\EE$, an 
externally generated electric
field are given. Denote by $\eta$  the
viscosity, $\rho$  the density, $\sigma^{-1}$  the electrical
resistivity and $\mu$ the magnetic permeability of the fluid (some positive constants).

To describe the interaction of the magnetic field and the
electrically conducting fluid in $\Omega$ mathematically  we combine the equations of 
the fluid
dynamics and electromagnetic field equations. We use Navier--Stokes equations to
model the fluid flow and include Lorenz force $\JJ \times \BB$ to express the
influence
 of the flow of  electric current across magnetic lines  on the fluid
motion
\begin{align}
\label{momentum.eq}
 - \eta \Delta{\uu} + \rho(\uu \cdot \nabla)\uu+\nabla{p} - \JJ \times \BB=
 {\bf F}.
\end{align}
The flow
of the electrically conducting fluid
across magnetic lines causes an electric current in the fluid. The electric
current alters the electromagnetic state of the system modifying the total
magnetic field, which creates the current in the fluid.
This phenomenon is expressed by means of Ohm's law
\begin{align}
\label{Ohm's.law} \sigma^{-1}\JJ+\nabla{\phi} - \uu \times \BB =\EE.
\end{align}
The total magnetic field $\BB=\BB_0+\cB(\JJ)$ is decomposed into a
sum of a given externally generated magnetic field $\BB_0$ and the
induced  magnetic field $\cB(\JJ)$, which is induced in the fluid  by the
electric current $\JJ$ caused by $\BB$. $\cB(\JJ)$
is a unique solution (see Lemma 2.2. in \cite{MS1}) of the
quasi-stationary form of Maxwell's equations
$$
 \nabla \times \mu^{-1} \cB(\JJ)=\JJ \ \ \
 \text{ and } \ \nabla \cdot \cB(\JJ)=0.
$$
It is also given by the Biot-Savart law
\begin{align}
\label{Biot-Savart}
 \cB(\JJ)({\bf x})=-\nabla \times \left(\frac{\mu}{4\pi}
 \int_{\Omega}\frac{ \JJ({\bf y})}{|{\bf x}-{\bf y}|}d {\bf y} \right)
 =-\frac{\mu}{4\pi}\int_{\Omega}\frac{{\bf x}-{\bf y}}{|{\bf x}-{\bf y}|^3}
 \times \JJ({\bf y}) d {\bf y},
\end{align}
for ${\bf x} \in \RR^3$. 
Note that using \eqref{Biot-Savart}, we are able to eliminate $\BB$ from
\eqref{momentum.eq}--\eqref{Ohm's.law} and solve for $\JJ$ instead.
Solving for $\BB$ we would be facing the problem that the magnetic
field is defined on all space and satisfies different equations
inside and outside the region $\Omega$: The Navier-Stoke's
equations are posed inside the region occupied by the fluid,
whereas the Maxwell's equations have to be solved in all of space.
Thus, the boundary conditions for $\BB$ on the surface of the
region must be specified, which is possible only for perfect
conductors or perfect insulators. This makes the prescribed
boundary conditions somewhat artificial.

We also consider the continuity equations
\begin{align}
\label{continuity.eq}
 \nabla \cdot \uu=0, \hskip5cm \nabla \cdot \JJ=0
\end{align}
to describe the incompressibility of the fluid and preservation of
charge.

Denote by $\nn$ the outward unit normal vector field on the
boundary of $\Omega$ and the stress tensor by
$$
T (\uu,p):=-p{\mathcal I}+\eta\left(\nabla \uu +(\nabla \uu)^T \right)=:-p{\mathcal I}+\eta D(\uu).
$$
Let the boundary, denoted by $\Gamma$, of
the domain $\Omega$ consist of
several relatively open and pairwise disjoint components
$\Gamma_i$'s, $i=1,\dots, 4$, i.e.
$\Gamma=\bar{\Gamma}_1 \cup \bar{\Gamma}_2 \cup \bar{\Gamma}_3 \cup \bar{\Gamma}_4$.
Different boundary conditions are prescribed on each component.
The boundary conditions we consider are \\
Dirichlet type (prescribed velocity)
\begin{align} \label{Dirichlet velocity}
 \uu={\bf g}_1 \ \ \text{on}  \ \ \Gamma_1,
\end{align}
Neumann type (prescribed stress)
\begin{align} \label{Neuman velocity}
  T \nn={\bf h}_2 \ \ \text{on} \ \ \Gamma_2
\end{align}
and mixed type for velocity and stress
\begin{eqnarray} \label{mixed}
 \begin{array}{c}
 \uu \cdot \nn =g_3 \quad \text{and} \quad
     (T \nn-((\nn T) \cdot \nn)\nn)={\bf h}_3 \ \ \text{on} \ \ \Gamma_3,
     \\ \\
 \uu- (\uu \cdot \nn)\nn ={\bf g}_4 \quad  \text{and}  \quad
     ( T \nn) \cdot \nn =h_4 \ \ \text{on} \ \ \Gamma_4.
\end{array} \end{eqnarray}
These boundary conditions are helpful in
modeling the free boundary value problems, when dealing with the artificially
truncated computational domains and the boundary conditions on the artificial
boundaries.

Let also the boundary $\Gamma$ consist of two relatively open and pairwise disjoint components
 $\Sigma_1$ and $\Sigma_2$, and consider \\
Neumann type (prescribed electric current flux through the walls) boundary condition
\begin{align}
\label{current flux}
 \JJ \cdot \nn =j \ \ \text{on} \ \ \Sigma_1
\end{align}
and Dirichlet type (prescribed electric potential) boundary condition 
\begin{align}
\label{elec potential}
  \phi=k \ \ \text{on} \ \ \Sigma_2.
\end{align}
This helps us to model two different cases: the external
magnetic field is given and   $\Sigma_1$ is not electrically
insulating then $j \not =0$; the magnetic field is generated by the external
conductors embedded into $\Sigma_2$. Note that incorporating the electric potential
is useful for various control problems.

\subsection{Function spaces}
Denote by $H^{s}(\Omega)$ the Sobolev space of square
integrable functions $v$ on $\Omega$ with square integrable distributional
derivatives $D^{\alpha}v$ up to order $s$ with the norm
$$
||v||_{H^{s}(\Omega)}=\Bigl( \sum_{|\alpha|\leq s} ||D^{\alpha}v||^2_
{L_2(\Omega)}\Bigr) ^{1 \over 2}.
$$
For vector--valued functions  $\vv=(v_1,v_2,v_3)$ define
\begin{eqnarray}\nonumber
 \begin{array}{c}
 \bH^s(\Omega):=\{\vv :v_i \in H^s(\Omega), i=1 \dots 3\}, \\
\text{{\bf L}}_2(\Omega):=\{\vv :v_i \in L_2(\Omega), i=1 \dots 3\ \}.
\end{array}
\end{eqnarray}
Let $ H^{1/2}(\Gamma)$ denote the fractional order Sobolev space and its dual
$H^{1/2}(\Gamma)'$.

The minimal regularity assumptions on the given data are
\begin{eqnarray} \nonumber
\begin{array}{c}
    {\bf F} \in \bH^1(\Omega)' \text{, }
    \ \ \ \EE \in \LL_{2}(\Omega), \ \ \ \BB_0 \in \bH^1(\Omega), \\ \\
     {\bf g}_1 \in \bH^{1/2}(\Gamma_1), \ \ \
     g_3 \in H^{1/2}(\Gamma_3), \ \ \ {\bf g}_4 \in \bH^{1/2}(\Gamma_4)
     \ \ \text{with} \ \ {\bf g}_4 \cdot \nn=0 \ \ \text{on} \ \ \Gamma_4, \\ \\
      {\bf h}_2 \in H^{1/2}(\Gamma_2)', \ \ \ {\bf h}_3 \in \bH^{1/2}(\Gamma_3)'
     \ \ \ h_4 \in H^{1/2}(\Gamma_4)' \ \ \text{with} \ \
     {\bf h}_3 \cdot \nn=0 \ \ \text{on} \ \ \Gamma_3.
\end{array}\end{eqnarray}
Note that \eqref{continuity.eq} imposes the following compatibility conditions on the
boundary data
\begin{eqnarray} \label{comp.u}
 \int_{\Gamma_1} {\bf g}_1 \cdot \nn +\int_{\Gamma_3} g_3=0, \quad
 \text{if} \ \ \Gamma_2 \cup \Gamma_4 = \emptyset, \quad \text{and} \quad
 \int_{\Sigma_1} j =0, \quad \text{if} \quad \Sigma_2=\emptyset.
\end{eqnarray}

The solution of (\ref{momentum.eq})--(\ref{elec potential}) is not unique if $\Gamma_2 \cup \Gamma_4 =\emptyset$
 and/or $\Sigma_2=\emptyset$. This is due to the fact that (\ref{momentum.eq})--(\ref{elec potential}) then
only involve the derivatives of $p$ and/or $\phi$. To avoid the
non-uniqueness  of the solution we solve for $p \in \dot{L}_2(\Omega)$ (subspace of
$L_2(\Omega)$ consisting of functions with mean zero)
and $\phi \in \dot{H}^1(\Omega)$ (subspace of
$H^1(\Omega)$ consisting of functions with mean zero). Therefore, we seek the solution
\begin{eqnarray} \nonumber
\begin{array}{c}
   \uu \in \bH^1(\Omega) \ \
       \text{ satisfying (\ref{momentum.eq})-(\ref{mixed})}, \\ \\
   \JJ \in \LL_2 \ \ \text{ satisfying } \ \
        \JJ \cdot \nn=j \ \ \text{on} \ \ \Sigma_1, \\ \\
  p\in M_p:=\left\{ \begin{array}{ccc} \dot{L}_2(\Omega)  &,& \text{if}
     \ \ \Gamma_2 \cup \Gamma_4 =\emptyset, \\
 L_2(\Omega) &,& \text{otherwise} \end{array} \right. \\ \\
 \phi \in M_\phi:=\left\{ \begin{array}{ccc} \dot{H}^1(\Omega)  &,& \text{if}
     \ \ \Sigma_2 =\emptyset, \\
\{ \phi \in H^1(\Omega): \phi \  \text{satisfying}\ \ \eqref{elec potential} \} &,& \text{otherwise}. \end{array} \right.
\end{array}\end{eqnarray}

To derive an equivalent to (\ref{momentum.eq})--(\ref{elec potential}) variational formulation, 
we choose
the test functions $\vv$ in the space
$$
  \bH^1_\Gamma(\Omega):=\{\vv \in \bH^1(\Omega): \vv|_{\Gamma_1} =0, \
   \vv \cdot \nn|_{\Gamma_3}=0, \ \vv-(\vv \cdot \nn)\nn|_{\Gamma_4}=0\},
$$
$\KK \in \LL_2(\Omega)$, $q \in M_q:=M_p$
and
\begin{eqnarray} \nonumber
\begin{array}{c}
 \psi \in M_\psi:=\left\{ \begin{array}{ccc} \dot{H}^1(\Omega)  &,& \text{if}
     \ \ \Sigma_2 =\emptyset, \\
 \{q \in H^1(\Omega): q|_{\Sigma_2} =0 \} &,& \text{otherwise.} \end{array} \right.
\end{array}\end{eqnarray}

To simplify the notation we introduce the product spaces
\begin{align*}
  \XX_{(\uu,\JJ)}&:=\bH^1(\Omega) \times \LL_2(\Omega), \ \
  \XX_{(\vv,\KK)}:=\bH^1_{\Gamma}(\Omega) \times \LL_2(\Omega), \\
  M_{(p,\phi)}&:= M_p\times M_\phi  \ \ \text{and} \ \ \
  M_{(q,\psi)}:=M_q \times M_\psi .
\end{align*}

\section{Weak Formulation}

In this section we shortly recall how to derive a variational problem equivalent to
(\ref{momentum.eq})--(\ref{elec potential}), (see \cite{CMS} for details).

Define a bilinear form $a_{0}:\XX_{(\uu,\JJ)} \times \XX_{(\uu,\JJ)}
\rightarrow \RR$, a trilinear form $a_{1}:\XX_{(\uu,\JJ)} \times\XX_{(\uu,\JJ)}
 \times \XX_{(\uu,\JJ)}
\rightarrow \RR$, a bilinear form $b: \XX_{(\uu,\JJ)} \times \left(L_2(\Omega) \times H^1(\Omega) \right) \rightarrow
\RR$ and a linear form $\chi: \XX_{(\uu,\JJ)} \rightarrow \RR$ by
\begin{align*}
a_{0}\bigl((\vv_{1},\KK_{1}),(\vv_{2},\KK_{2})\bigr)&:= \frac{\eta}{2}
\int_{\Omega} D(\vv_1):D(\vv_2) + \sigma^{-1} \int_{\Omega}\KK_1 \cdot \KK_2 \cr
&+\int_{\Omega}\bigl((\KK_2\times \BB_0)\cdot \vv_1-(\KK_1\times \BB_0)\cdot
\vv_2\bigr),\cr
a_{1}\bigl((\vv_1,\KK_1),(\vv_2,\KK_2),(\vv_3,\KK_3)\bigr)&:=
\rho \int_{\Omega} \bigl((\vv_1\cdot \nabla)\vv_2\bigr) \cdot \vv_3 \cr
&+\int_{\Omega} \Bigl(\bigl(\KK_3 \times \cB(\KK_1)\bigr)\cdot \vv_2-
\bigl(\KK_2\times \cB(\KK_1)\bigr)\cdot \vv_3 \Bigr) ,\cr
b\bigl((\vv,\KK)(q,\psi)\bigr)&:=-\int_{\Omega}(\nabla \cdot \vv)q+
\int_{\Omega} \KK\cdot (\nabla \psi)
\end{align*}
and
$$
\chi (\vv,\KK)=\int_{\Omega} {\bf F} \cdot \vv + \int_{\Omega} \EE\cdot \KK +
\int_{\Gamma_2} {\bf h}_2\cdot \vv  +
\int_{\Gamma_3} {\bf h }_3\cdot \vv  +
\int_{\Gamma_4} h_4 \nn \cdot \vv .
$$

Let $\vv \in \bH_{\Gamma}^{1}(\Omega)$ and $\KK \in \LL_{2}(\Omega)$
be arbitrary test functions. Multiplying (\ref{momentum.eq}),
(\ref{Ohm's.law}) and \eqref{continuity.eq} by corresponding
test functions, integrating by parts and using some of the boundary
conditions on $\uu$ and $\JJ$ we get the following equivalent variational problem.

\vspace{0.3cm}\noindent {\bf Problem 1:} {\it Find $(\uu,\JJ) \in  \XX_{(\uu, \JJ)}$
s.th. $\uu|_{\Gamma_1}=\g_1, \ \ (\uu\cdot \nn)|_{\Gamma_3}=g_3, \ \
\bigl(\uu-(\uu\cdot \nn)\nn\bigr)|_{\Gamma_4}=\g_4$
and $(p,\phi) \in M_{(p,\phi)}$ such that
$$
a_{0}\bigl((\uu,\JJ),(\vv,\KK)\bigr)+
a_{1}\bigl((\uu,\JJ),(\uu,\JJ),(\vv,\KK)\bigr)
+b\bigl((\vv,\KK),(p,\phi)\bigr)=\chi (\vv,\KK),
$$
and
$$
b\bigl((\uu,\JJ),(q,\psi)\bigr)=\int_{\Sigma_1}j\psi
$$
for all $(\vv,\KK) \in  \XX_{(\vv,\KK)}$ and $(q,\psi) \in  M_{(q,\psi)}$.} \\

Next, we list some properties (proved in \cite[Lemma 2.2]{MS1}
and \cite[Corollary 1]{CMS}) of  $a_0$, $a_1$ and $b$ .

\begin{lemma}
\label{lemma:propertiesforms} \item a) The forms $a_0$, $a_1$, $b$ are bounded on
$\XX_{(\uu,\JJ)} \times \XX_{(\uu,\JJ)}$
, $\XX_{(\uu,\JJ)} \times \XX_{(\uu,\JJ)} \times \XX_{(\uu,\JJ)}$
and $\XX_{(\uu,\JJ)} \times M_{(p,\phi)}$, respectively, with
\begin{eqnarray*}
 \|a_1\| \lesssim \max\{\rho, \mu\}.
\end{eqnarray*}
\item b) If the intersection of $\bH^1_\Gamma(\Omega)$ and the null space $N(D)$ of the deformation
tensor $D$ is $\{0\}$, then the form $a_0$ is positive definite on $\XX_{(\vv,\KK)} \times \XX_{(\vv,\KK)}$, i.e.
$$
 a_{0}\bigl((\vv,\KK),(\vv,\KK)\bigr) \ge \alpha_0 \|(\vv,\KK)\|_{\XX_{(\vv,\KK)}}
$$
with $\alpha_0:=c(\Omega)\min\{\eta, \sigma^{-1}\}$, $c(\Omega)$ some
positive constant depending only on the domain.  
\item c) The bilinear form $b$ satisfies the {\it inf-sup} condition
$$
 \inf_{(q,\psi) \in M_{(q,\psi)}} \sup_{(\vv,\KK) \in \XX_{(\vv, \KK)}}
 \frac{b\bigl((\vv,\KK),(q,\psi)\bigr)}{\|(\vv,\KK)\|_{\XX_{(\vv, \KK)}}
 \|(q,\psi)\|_{M_{(q,\psi)}}} > 0.
$$
\end{lemma}
Note that $\bH^1_\Gamma(\Omega) \cap N(D)=\{0\}$, for example, if $\Gamma_1 \not= \emptyset$. 
Otherwise, the subsequent analysis would require the introduction of a suitable quotient space
of $\bH^1_\Gamma(\Omega)$.

Due to \cite[Corollary 1]{CMS} there exist the liftings of the boundary data
\begin{align*}
 & \uu_0 \in \bH^1(\Omega) \ \ \text{with} \ \nabla \cdot \uu_0=0 \
 \text{and} \ \uu_0|_{\Gamma_1}={\bf g}_1, \uu_0 \cdot \nn|_{\Gamma_3}=g_3
 \ \text{and} \ \uu_0-(\uu_0 \cdot \nn)\nn|_{\Gamma_4}={\bf g}_4, \cr
 & \JJ_0 \in \LL_2(\Omega) \ \ \text{with} \ \nabla \cdot \JJ_0=0 \
 \text{and } \ \ \JJ_0 \cdot \nn=j \ \text{on} \ \Sigma_1.
\end{align*}
Define $\hat{\uu}:=\uu-\uu_0$, $\hat{\boldsymbol{J}}:=\JJ-\JJ_0$. Also set
$\hat{\phi}:=\phi - \phi_0$ and $ \hat{p} := p -  p_0$, where
 $p_0=0$,
if $\Gamma_2 \cup \Gamma_4 \not =\emptyset$, and $\displaystyle{p_0=\frac{1}{|\Omega|}
\int_{\Omega}p}$,
otherwise, and $\phi_0$ is the $H^1-$lifting of the boundary data $k$ on $\Sigma_2$,
if $\Sigma_2 \not =\emptyset$, and $\displaystyle{\phi_0=\frac{1}{|\Omega|}\int_{\Omega}\phi}$,
otherwise. Define the bilinear form $a: \XX_{(\vv,\KK)} \times \XX_{(\vv,\KK)} \rightarrow \RR$ and the linear
form $\ell:\XX_{(\vv,\KK)}  \rightarrow \RR$ by
\begin{eqnarray} \nonumber
\begin{array}{c}
 a \bigl((\vv_1, \KK_1),(\vv_2,\KK_2)\bigr):=
 a_0\bigl((\vv_1, \KK_1),(\vv_2,\KK_2)\bigr)
 +a_1\bigl((\vv_1, \KK_1),(\uu_0,\JJ_0),(\vv_2,\KK_2)\bigr)\\
 +a_1\bigl((\uu_0,\JJ_0),(\vv_1, \KK_1),(\vv_2,\KK_2)\bigr) \\
\end{array} \end{eqnarray}
and
\begin{eqnarray} \nonumber
\begin{array}{c}
 \ell(\vv,\KK):=\chi (\vv,\KK) -a_0\bigl((\uu_0,\JJ_0),(\vv,\KK)\bigr)-
          a_1\bigl((\uu_0,\JJ_0),(\uu_0,\JJ_0),(\vv,\KK)\bigr) 
          -b((\vv,\KK),(p_0,\phi_0)).
\end{array} \end{eqnarray}
Note that if $(\uu_0,\JJ_0)$ have small norm in $\XX_{(\uu,\JJ)}$, 
then the form $a$ {\bf is coercive on} $\XX_{(\vv,\KK)}$, i.e.
\begin{equation} \label{coercive}
 a((\vv,\KK),(\vv,\KK)) \ge \Big(\alpha_0-2\|a_1\| \ \|(\uu_0, \JJ_0)\|_{\XX_{(\uu,\JJ)}}\Big)
 \|(\vv,\KK)\|^2_{\XX_{(\vv,\KK)}}
\end{equation}
for all $(\vv,\KK) \in \XX_{(\vv,\KK)}$. Substituting $\uu=\hat{\uu}+\uu_0$, $ \JJ =\hat{\JJ}+\JJ_0$,
$ \phi =\hat{\phi} + \phi_0$ and $ p = \hat{p}+p_0$ into {\bf Problem 1} we get its
equivalent formulation.

\vspace{0.3cm} \noindent {\bf Problem 2:} {\it Find $(\hat{\uu}, \hat{\JJ}) \in
\XX_{(\vv,\KK)}$ and $(\hat{p},\hat{\phi}) \in M_{(q,\psi)}$
satisfying
$$
a \bigl((\hat{\uu},\hat{\JJ}),(\vv,\KK)\bigr)+
a_1\bigl((\hat{\uu},\hat{\JJ}),(\hat{\uu},\hat{\JJ}),(\vv,\KK)\bigr)+
b\bigl((\vv,\KK), (\hat{p},\hat{\phi})\bigr)=\ell(\vv,\KK)
$$
for all $(\vv,\KK) \in \XX_{(\vv,\KK)}$ and
$ \quad
b\bigl((\hat{\uu},\hat{\JJ}),(q,\psi)\bigr)=0 \ \ \text{for all} \ \ (q,\psi)
\in M_{(q,\psi)}.
$} \\

Next, define
$$
 \bV:=\left\{(\vv,\KK) \in \XX_{(\vv,\KK)}: b\bigl((\vv,\KK),(q,\psi)\bigr)=0 \text{
for all} \ (q,\psi) \in M_{(q,\psi)} \right\} $$ and consider

\vspace{0.3cm} \noindent {\bf Problem 3:} {\it Find $(\hat{\uu},\hat{\JJ}) \in \bV$ such that
$$
 a \bigl((\hat{\uu},\hat{\JJ}),(\vv,\KK)\bigr)+
 a_1\bigl((\hat{\uu},\hat{\JJ}),(\hat{\uu},\hat{\JJ}),(\vv,\KK)\bigr)
 =\ell(\vv,\KK)
$$
for all $(\vv,\KK) \in \bV$.} \\
The standard nonlinear version of the classical LBB
(Ladyzhenskaya-Babushka-Brezzi) theory (see Chapter 4.1 in
\cite{GR}) tells us that {\bf Problem 2} and {\bf Problem 3} are
equivalent, if the form $b$ satisfies the {\it inf-sup} condition,
and that {\bf Problem 3} is uniquely solvable, if the form $a$ is
coercive and bounded. In other words, due to Lemma
\ref{lemma:propertiesforms} and \eqref{coercive}, the LBB-theory
allows us to further transform {\bf Problem 2} and solve {\bf
Problem 3} for the unknown velocity $\hat{\uu}$ and electric
current density $\hat{\JJ}$ on a subspace $\bV$ of
$\XX_{(\vv,\KK)}$. Thus, if $(\hat{\uu},\hat{\JJ}) \in \bV$ is a
solution of {\bf Problem 3}, then there exist a unique pair
$(\hat{p},\hat{\phi}) \in M_{(q,\psi)}$ such that
$(\hat{\uu},\hat{\JJ}, \hat{p},\hat{\phi})$
 solves {\bf Problem 2}. And, given
any $((\hat{\uu},\hat{\JJ}),(\hat{p},\hat{\phi})) \in
\XX_{(\vv,\KK)} \times M_{(q,\psi)}$, solution of  {\bf Problem 2}, then
$(\hat{\uu},\hat{\JJ})$ is in $\bV$ and  solves {\bf Problem 3}. By
\cite[Theorem in Section 3]{CMS} {\bf Problem 3} is uniquely solvable if
\begin{equation}\label{bound:u0J0}
  \|(\uu_0,\JJ_0)\|_{\XX_{(\uu,\JJ)}}\ < \frac{\alpha_0}{2\|a_1\|}
\end{equation}
and
\begin{equation}\label{bound:ell}
 \|\ell\|_{\bV'} < \frac{(\alpha_0-2\|a_1\|
 \|(\uu_0,\JJ_0)\|_{\XX_{(\uu,\JJ)}})^2}{4\|a_1\|}.
\end{equation}
The definition of $\alpha_0=c(\Omega)\min\{\eta, \sigma^{-1}\}$
implies that any given set of data will satisfy
\eqref{bound:u0J0}-\eqref{bound:ell} if the viscosity $\eta$ and
the electric resistivity $\sigma^{-1}$ are large enough.

\section{From weak formulation to  frame discretization}

We start by reformulating {\bf Problem 3} to fit the following abstract setting: There 
exists a separable Hilbert space $\mathcal{H}$ such that $\mathbf{V} \subset \cH
\subset \mathbf{V}'$ with bounded and dense inclusions. The triple
$(\mathbf{V}, \cH, \mathbf{V}')$ is then called a Gelfand triple.
The duality between $\bV'$ and $\bV$ is identified on $\cH$ using
the inner product $\langle \cdot,\cdot \rangle_{\cH}$ of $\cH$. 
There exists an operator $A:\bV \rightarrow \bV'$ such
that $\langle A v,w\rangle_{\bV' \times \bV}:=a(v,w)$ defines an elliptic bilinear
form, i.e., there exist positive constants $\alpha$, $\beta$ such
that $\alpha \|v\|_{\bV}^2 \le a(v,v) \le \beta \|v\|_{\bV}^2$ for all
$v \in \bV$.
The ellipticity of $a$ implies that $\|A v\|_{\bV'} \sim
\|v\|_{\bV}$ and that $A$ is a boundedly invertible operator with 
$\|A^{-1}\| \le \alpha^{-1}$. We also assume that there
exists a trilinear form $a_1$ inducing a bounded bilinear operator
$A_1:\bV\times\bV \rightarrow \bV'$, defined by $\langle
A_1(v,w),z\rangle_{\bV' \times \bV} := a_1(v,w,z)$ for $v,w,z \in \bV$. It
has been shown in \cite{CMS} that {\bf  Problem 3}  translates into the
following abstract problem:

\vspace{0.3cm}\noindent \vspace{0.3cm}{\it  Find
$u:=(\hat{\uu}, \hat{\JJ}) \in \bV$ such that
\begin{equation}
\label{abfor}
A u +A_1(u,u)= \ell,
\end{equation}
where $\ell \in \bV'$ is a functional on ${\mathbf V}$.}
\vspace{0.3cm}

The solvability of \eqref{abfor} is ensured by \cite[Lemma 2]{CMS} and summarized by the following theorem.
\begin{tm}
\label{tmfixpt}
Let $H: \bV \rightarrow \bV$ be given by  $H(v):=A^{-1}\L \ell - A_1(v,v)\R$. The following hold
\begin{itemize}
\item[(a)] If $0< r < \frac{\alpha}{2 \|A_1\|}$, then $H_{| B_r}$ is a contraction with 
Lipschitz constant $L:=\frac{2 \|A_1\| r}{\alpha}$,
where $B_r \subset \bV$ is the closed ball of radius $r$ centered at the origin;
\item[(b)] If $0< r < \frac{\alpha}{\|A_1\|}$ and $\|\ell\|_{\bV'} \leq r(\alpha -\|A_1\| r)$, then $H_{|B_r} :B_r \rightarrow B_r$;
\item[(c)] If $\|\ell\|_{\bV'} \leq \frac{\alpha^2}{4 \|A_1\|}$, then   \eqref{abfor} has 
a unique solution $u$ with
$\|u\|_\bV < \frac{\alpha}{2 \|A_1\|}$ and the solution is given
by the fixed point iteration
\begin{eqnarray}
\label{fixpt}
u_{n+1} &=& H u_{n},  \quad u_0 = 0, \quad n \in \NN_0, \\
u&=& \lim_{n \rightarrow \infty} u_n. \notag
\end{eqnarray}
\end{itemize}
\end{tm}

Note that \eqref{fixpt} is equivalent to
\begin{equation}
\label{ellpt}
A u_{n+1} = \ell -A_1(u_n,u_n), \quad n \in \NN_0, \quad u_0=0.
\end{equation}
Under our assumptions on $A$, the equations in (\ref{ellpt}) are elliptic operator equations.
We show in this section how to derive a discrete problem  equivalent to the abstract nonlinear problem in \eqref{abfor} and prove a result
similar to Theorem \ref{tmfixpt} showing that the
solution of the discrete problem exists and is unique under certain assumptions on the 
parameters of the
original MHD problem.
The discrete problem is obtained using suitable {\it stable}, {\it redundant}, and 
{\it nonorthogonal} expansions, so--called Gelfand frames 
for the Gelfand triple $(\mathbf{V}, \cH, \mathbf{V}')$. We also show that the corresponding 
discrete fixed point iteration can be numerically realized  efficiently. In 
particular, we present
the realization of the key numerical routine {\bf SOLVE} used in our discrete fixed point 
iteration to approximate  adaptively
the solution of the elliptic problems in (\ref{ellpt}).

\subsection{Gelfand Frames}
\label{prelim}

In the following, the sequence space  $\ell_2(\mathcal{N})$ on the 
countable index set ${\mathcal N} \subset \RR^d$ is
induced by the norm
$$
 \|\vec{\mathbf{c}}\|_{\ell_2(\cN)}:=\left(\sum_{n\in\mathcal{N}}|c_n|^2 \right)^{1/2}, \quad \vec{\mathbf{c}}=\{c_n\}_{n \in \cN} \in 
 \ell_2(\cN).
$$
The space $\ell_0(\mathcal{N}) \subset \ell_2(\mathcal{N})$ is the subspace of sequences with compact support.
Denote $\langle\cdot,\cdot\rangle_\cH$ and $\|\cdot\|_{\mathcal H}$ the inner product and 
the norm on
the separable Hilbert space $\mathcal{H}$, respectively. A sequence $\mathcal{F}:=\{f_n\}_{n \in \mathcal{N}}$ in $\mathcal{H}$ is a \textit{frame} for $\mathcal{H}$ if
\begin{equation}
\label{framestability}
\|f\|_{\mathcal H}^2 \sim \sum_{n \in \cN} \bigl|\langle f, f_n \rangle_{\cH}\bigr|^2, \quad \text{for all }f\in \mathcal{H}.
\end{equation}
Due to \eqref{framestability}  the corresponding operators of analysis and synthesis
given by
\begin{equation}
\label{analysisop}
F:\mathcal{H}\rightarrow\ell_2(\mathcal{N}),\quad f\mapsto\bigl(\langle f,f_n\rangle_{\cH}\bigr)_{n\in\cN},
\end{equation}
\begin{equation}
\label{synthesisop}
F^*:\ell_2(\mathcal{N})\rightarrow\mathcal{H},\quad \vec{\mathbf{c}}\mapsto\sum_{n\in\cN} c_nf_n,
\end{equation}
are bounded.
The composition $S:=F^*F$ is a boundedly invertible (positive and self--adjoint) operator 
called the {\it frame operator} and $\tilde{\mathcal{F}}:=\{S^{-1} f_n\}_{n  \in \cN}$ 
is again a frame for $\mathcal{H}$, called the {\it canonical dual frame}, with corresponding analysis and synthesis operators
\begin{equation}
\label{dualframeops}
\tilde F: = F(F^* F)^{-1}, \quad \tilde F^*:=(F^* F)^{-1} F^*.
\end{equation}
In particular, one has the following orthogonal decomposition of $\ell_2(\mathcal{N})$
\begin{equation}
\label{decomp}
\ell_2(\mathcal{N})=\Ran(F)\oplus\Ker(F^*),
\end{equation}
and
\begin{equation}
\label{Qdef}
\mathbf{Q} := F(F^* F)^{-1} F^* :\quad\ell_2(\mathcal{N}) \rightarrow\Ran(F),
\end{equation}
is the orthogonal projection onto $\Ran(F)$.

The frame $\mathcal{F}$ is a Riesz basis for
$\mathcal{H}$ if and only if $\Ker(F^*)=\{0\}$. In general, we
assume that $\{0\}$ is a proper subspace of $\Ker(F^*)$. In other words, due to the 
redundancy of the frame there may exist sequences 
$\vec{\mathbf{c}}=\{c_n\}_{n \in \cN} \neq \vec{\mathbf d}=\{d_n\}_{n \in \cN}$ in 
$\ell_2(\cN)$
such that $\displaystyle{\sum_{n \in \cN} c_n f_n = \sum_{n \in \cN} d_n f_n}$. 
In particular, the redundancy may lead to the situation when
 a small perturbation  $ \vec{\mathbf{d}}$ of
the coefficient sequence  $\vec{\mathbf{c}}$ has no effect on the
synthesis operator. This possible reduction effect on  errors, noise, and numerical
round-offs is the motivation for using frames for the
applications, where tolerance to errors is required.
 The intrinsic stability of frames is also expected to play a role in the conditioning of the discretizations
of operator equations and leads to additional robustness that the discretization inherits from the frame.
In fact, it has been recently confirmed by numerical experiments in \cite{W} that increasing (uniform) redundancy of
the frame one improves  the conditioning of the corresponding discretization matrices.

It is somewhat perplexing that different sets of coefficients yield equivalent representations
 of the {\it same} element of $\cH$.
It is not clear then what
 are the ``good and computable coefficients'': The importance of the canonical dual frame is its use 
 in reconstruction of
any   $f \in \cH$, i.e. 
\begin{equation}
\label{reconstruction}
f =S S^{-1} f = \sum_{n\in\cN} \langle f, S^{-1} f_n \rangle_\cH f_n=S^{-1} S f = \sum_{n\in\cN} \langle f, f_n \rangle_\cH S^{-1} f_n.
\end{equation}
Since a frame is typically overcomplete in the sense that the coefficient functionals 
$\vec{\mathbf{c}}=\{c_n(f)\}_{n \in \mathcal{N}} \in \ell_2(\cN)$ in the representation
\begin{equation}
\label{coefffunctionals}
f = \sum_{n \in \mathcal{N}} c_n(f) f_n
\end{equation}
are in general not unique ($\Ker(F^*)\neq\{0\}$), there exist many possible non--canonical duals $\{\tilde f_n\}_{n\in\mathcal{N}}$ in $\mathcal{H}$ for which
\begin{equation}
\label{reconstnoncanonical}
f = \sum_{n \in \mathcal{N}} \langle f, \tilde f_n \rangle_\cH f_n.
\end{equation}
\\

A more general definition of frames is required for Banach spaces.
For details on Banach frames we refer, for example, to \cite{DFR,FG3,FG,G0,G1}.
Throughout this paper we  make use of {\it Gelfand frames} that are particular instances of 
Banach frames. So we do not introduce the latter here in full generality.
Assuming that $\cB$ is a Banach space continuously and densely embedded in $\mathcal{H}$ we get
\begin{equation}
\label{gelfandtriple}
\cB \subseteq \mathcal{H}\simeq\mathcal{H}' \subseteq \cB'.
\end{equation}
If the right inclusion is dense, then $(\cB,\cH,\cB')$ is called a \textit{Gelfand triple}.
The symbol $\simeq$ stands for the canonical Riesz identification of $\mathcal{H}$ with its dual $\mathcal{H}'$.
\begin{definition} \label{gelfrm} A frame $\cF$ (here $\tilde \cF$ is some dual frame, e.g., the canonical dual frame) for $\cH$ is called a \textit{Gelfand frame} for the Gelfand triple $(\cB,\cH,\cB')$, if $\cF\subset\cB$, $\tilde \cF \subset \cB'$ and there exists a Gelfand triple $\bigl(\cB_d,\ell_2(\cN),\cB_d'\bigr)$ of sequence spaces such that
\begin{equation}
\label{gelfandframeprop}
F^*:\cB_d \rightarrow \cB,\; F^* \vec{\mathbf{c}}=\sum_{n\in\cN} c_n f_n \quad\text{and}\quad \tilde F:\cB \rightarrow \cB_d,\; \tilde Ff=\bigl(\langle f,\tilde f_n\rangle_{\cB\times\cB'}\bigr)_{n\in\cN}
\end{equation}
are bounded operators.\\
\end{definition}

\begin{rems}

1. If $\cF$ (again  $\tilde \cF$ is some dual frame, e.g., the canonical dual frame) is a Gelfand frame for the Gelfand triple $(\cB,\cH,\cB')$
with respect to the Gelfand triple of sequences $\bigl(\cB_d,\ell_2(\cN),\cB_d'\bigr)$, then by duality also the operators
\begin{equation}
\label{duality}
\tilde F^*:\cB_d' \rightarrow \cB',\;\tilde F^*\vec{\mathbf{c}}=\sum_{n\in\cN} c_n \tilde f_n \quad\text{and}\quad F:\cB'
\rightarrow \cB_d',\;Ff=\bigl( \langle f,f_n\rangle_{\cB'\times\cB} \bigr)_{n \in \mathcal{N}}
\end{equation}
are bounded, see, e.g., \cite{hackbusch} for details.

2. If $\cB  =\mathcal H $ then Definition \ref{gelfrm} becomes the definition of  frames for  Hilbert spaces.

3. Even if the solution space $\mathbf{V} \subset \cH= \mathbf{L}^2(\Omega) \subset \mathbf{V}'$ is a 
Hilbert space, by using the notation ``$\cB$'' here we want to emphasize that the frame $\cF$ we consider 
is {\it not} a Hilbert space frame for $\mathbf V$. 
It is a frame for $\cH$, which also characterizes $\mathbf V$ (as a subspace of $\cH$) with 
 frame coefficients $\langle f, \tilde f_i \rangle_{\cH}$ belonging to some suitable sequence space 
 $\mathbf V_d \subset \ell_2(\cN)$ and being computed using $\cH$-inner product. Of course, we can 
 consider genuine Hilbert frame expansions for $\mathbf{V}$ (as it is done, for example, in 
 \cite{S}). This may create some confusion because one might be tempted to use the  $\bV-$inner product when determining the duality pair. 
This, in most of the cases,  is not compatible with  numerical implementations.

4. Definition \ref{gelfrm} generalizes the following case to pure frames:
Consider a  wavelet system
$\Psi:=\{\psi_{j,k}\}_{j \geq -1, k \in \mathcal{J}_j}$ on 
$\Omega$ ($\mathcal{J}_j$ is a suitable set of indexes depending on the scale 
$j$, see \cite{U} for details),  
$\cB = H^s(\Omega)$, $\cH = L_2(\Omega)$ and 
$$\cB_d = \ell_{2,2^{s \cdot}} := 
\{\vec{\mathbf{d}} := \{d_{j,k}\}_{j \geq -1, k  \in \mathcal{J}_j}: \left 
(\sum_{j \geq -1} \sum_{k \in \mathcal{J}_j} 2^{2 s j} |d_{j,k}|^2 \right )^{1/2}<\infty \}.$$
It is well known that if $\Psi$ is a Riesz basis for $L_2(\Omega)$ and its elements, together with those of its biorthogonal dual basis $\tilde \Psi:=\{\tilde \psi_{j,k}\}_{j \geq -1, k \in \mathcal{J}_j}$, are compactly supported, smooth enough, and with a sufficient number of vanishing moments, then $H^s(\Omega)$ is fully characterized by $\Psi$ in the sense that $f \in H^s(\Omega)$ if and only if
\begin{equation}
 f  = \sum_{j \geq -1} \sum_{k \in \mathcal{J}_j} \langle f, \tilde \psi_{j,k} \rangle_{L_2(\Omega)} \psi_{j,k}
\end{equation}
and
\begin{equation}
 \|f\|_{H^s(\Omega)} \sim  \left (\sum_{j \geq -1} \sum_{k \in \mathcal{J}_j} 2^{2 s j} |\langle f, \tilde \psi_{j,k} \rangle_{L_2(\Omega)}|^2 \right )^{1/2}.
\end{equation}
See \cite{DFR} for the same characterization by using pure wavelet frames, constructed by Overlapping Domain 
Decomposition. Note that there exists a natural unitary isomorphism from $\ell_{2,2^{s \cdot}}$ into $\ell_2$ 
given by
\begin{equation}
D_{H^s(\Omega)} :\ell_{2,2^{j s \cdot}} \rightarrow \ell_2, \quad \vec{\mathbf{d}} := 
\{d_{j,k}\}_{j \geq -1, k \in \mathcal{J}_j} \mapsto D_{H^s(\Omega)} \vec{\mathbf{d}} := 
\{2^{j s} d_{j,k}\}_{j \geq -1, k \in \mathcal{J}_j}. 
\end{equation}
\end{rems}

Keeping in mind the example in the above REMARK 4., we  proceed to the numerical treatment of abstract 
elliptic operator equations by means of Gelfand frame discretizations. 

\subsection{Adaptive Numerical Frame Schemes for Elliptic Operator Equations} \label{resolution}
To implement the fixed point iteration described in Theorem
\ref{tmfixpt}, we first study the solvability (for fixed
$u^{(n)}$) of the linear  operator equations in \eqref{ellpt}.
Generally, such equations are of the form
\begin{equation} \label{basicproblem}
A u=f,
\end{equation}
where $A$, as before, is  a boundedly invertible operator
from Hilbert space $\bV$ into its dual $\bV'$,
\begin{equation} \label{2.2}
  \norm{Au}_{\bV'} \sim \norm{u}_{\bV}, \quad u \in \bV.
\end{equation}
We also have that
\begin{equation} \label{bilin}
a(v,w):=\langle {A} v, w\rangle_{\bV' \times \bV} \ ,
\end{equation}
defines a bilinear form on $\bV$, where $\langle\cdot,\cdot\rangle_{\bV \times \bV'}$ defines the 
dual pairing of $\bV$ and $\bV'$. The
form $a$ is {\em elliptic},  i.e., there exist positive constants $\alpha$, $\beta$ such
that 
\begin{equation}
\label{elliptic}
\alpha \|v\|_{\bV}^2 \le a(v,v) \le \beta \|v\|_{\bV}^2
\end{equation} 
for all
$v \in \bV$, and $a$ is {\em non-symmetric}. The assumption on the non-symmetry  of $a$ is motivated
by the MHD example presented in Sections 2-3.

Here and throughout the rest of the paper we assume  that $\cF=\{f_n\}_{n\in\cN}$ is a Gelfand frame 
for the Gelfand triple $(\bV, \cH,\bV')$ with
$(\bV_d,\ell_2(\cN),\bV_d')$ being the corresponding Gelfand triple of sequence spaces. Moreover, 
keeping in mind REMARK 4. in Subsection 4.1, we also assume that there exists a  unitary 
isomorphism $D_{\bV}: \bV_d \rightarrow \ell_2(\cN)$, so that its
$\ell_2(\cN)$--adjoint $D_{\bV}^*:\ell_2(\cN) \rightarrow {\bV}_d'$
is also an isomorphism.

We  show,  next, how the Gelfand frame setting can be used for the adaptive
numerical treatment of elliptic operator equations   \eqref{basicproblem}.
 Following, e.g. \cite{CDD1,DFR,S}, one uses  frame expansions to convert the problem \eqref{basicproblem}
into an operator equation on $\ell_2(\cN)$. The problem that arises is that the redundancy of the frame leads to 
a singular discretization matrix. Nevertheless, in Theorem \ref{uvectheorem} below we show that this can
be
handled in practice and that the solution of  \eqref{basicproblem} can be computed by a version 
of  
Richardson iteration applied to the associated normal equations. The resulting scheme is not 
directly implementable since one has
to deal with infinite matrices and vectors. Therefore, similarly to 
\cite{CDD1, CDD2, S}, we also show how the scheme can be
transformed into an implementable scheme using ``finite''
versions of the  building blocks procedures we introduce in Subsection 4.2.2. The result is a
convergent adaptive frame algorithm.

\subsubsection{A series representation}
We start by generalizing Lemma 4.1 and Theorem 4.2 given in \cite{DFR} to the case of non--symmetric $a$.
 We give the
detailed proofs to emphasize the difference between the symmetric
and non--symmetric cases.

\begin{lemma}\label{lemma:invetibilitybfA} Under the assumptions \eqref{bilin}, \eqref{elliptic} on $A$, the operator
\begin{equation}
\label{Gdef}
\ba := (D_{\bV}^*)^{-1} F A F^* D_{\bV}^{-1}
\end{equation}
is a bounded operator from $\ell_2(\cN)$ to $\ell_2(\cN)$.
Moreover $\ba$ is boundedly invertible on its range $\Ran(\ba) =
\Ran((D_{\bV}^*)^{-1} F)$.
\end{lemma}
\begin{proof}
Since $\ba$ is a composition of bounded operators
$D_{\bV}^{-1}:\ell_2(\cN) \rightarrow \bV_d$,
$F^* : \bV_d \rightarrow \bV$,
$A:\bV \rightarrow \bV'$,
$F:\bV' \rightarrow \bV_d'$ and
$(D_{\bV}^*)^{-1} : \bV_d' \rightarrow \ell_2(\cN)$,
$\ba$ is a bounded operator from $\ell_2(\cN)$ to $\ell_2(\cN)$. Moreover, from the decomposition \eqref{Gdef}  we get
\begin{equation}
\label{kerG}
\Ker(\ba) = \Ker(F^*D_{\bV}^{-1}), \quad \Ran(\ba) = \Ran((D_{\bV}^*)^{-1} F).
\end{equation}
Define ${\bf L} := (D_{\bV}^*)^{-1} F
F^* D_{\bV}^{-1}$. Note that $\Ker({\bf L}) =
\Ker(F^*D_{\bV}^{-1})$ and
 $\Ran({\bf L}) = \Ran((D_{\bV}^*)^{-1} F)$. The fact that $
 \ell_2(\cN)=\Ker({\bf L}^*)\oplus \Ran({\bf L})$
implies, due to the self-adjointness of ${\bf L}$, that
\begin{equation}
\label{Gdecomp}
\ell_2(\cN)=\Ker(F^*D_{\bV}^{-1})\oplus \Ran((D_{\bV}^*)^{-1} F).
\end{equation}
Therefore,
\begin{equation}
\label{Gspaces}
\ba_{|\Ran(\ba)}:\;\Ran(\ba) \rightarrow \Ran(\ba)
\end{equation}
is boundedly invertible. 
\end{proof}

Denote by  $ \mathbf{P}:\ell_2(\cN)\to\Ran(\ba)$ the orthogonal projection of $\ell_2(\cN)$ onto $\Ran(\ba)$.

\begin{tm}  \label{uvectheorem}
Let $A$ satisfy \eqref{bilin} and \eqref{elliptic}. Denote
\begin{equation}
\label{fvecdef}
\vec{\mathbf{f}}:=(D_{\bV}^*)^{-1} F f
\end{equation}
and $\ba$ as in \eqref{Gdef}. Then the solution $u$ of \eqref{basicproblem}
can be computed by
\begin{equation}
\label{ucompute}
u = F^* D_{\bV}^{-1} \mathbf{P} \vec {\mathbf u }
\end{equation}
where $\vec{\mathbf{u}}$ solves
\begin{equation}
\label{uvecproblem}
 \mathbf{P}  \vec{\mathbf{u}} = \L \alpha^* \sum_{n=0}^\infty (\id -\alpha^* \ba^* \ba)_{| \Ran(\ba)}^n \R \ba^* \vec{\mathbf{f}},
\end{equation}
with $0 < \alpha^* < 2/\lambda_{\max}$, where $\lambda_{max}=\|\ba^* \ba\|_2$  with $\| \cdot\|_2$ being the usual
spectral norm.
\end{tm}
\begin{proof}
We have $u=\sum_{n\in\cN} \langle u, \tilde f_n \rangle_\cH f_n$ in $\cH$.
Since $\cF$ is a Gelfand frame, $F^*\tilde F:\bV\to\bV$ is bounded
and implies $u = F^* \tilde F u=\sum_{n\in\cN} \langle u, \tilde f_n \rangle_{\bV\times\bV'} f_n$ in $\bV$.
Moreover,  \eqref{basicproblem} is equivalent to the following system of equations
\begin{equation}
\label{opeqdiscrete}
\sum_{n\in\cN} \langle u, \tilde f_n \rangle_{\bV\times\bV'}\langle A f_n,f_m\rangle_{\bV'\times\bV} = \langle f,f_m\rangle_{\bV'\times\bV}, \quad m\in\cN.
\end{equation}
Denote $\vec{\mathbf{u}} := D_{\bV} \tilde F u$ and $\vec{\mathbf{f}}$, $\ba$ as in \eqref{fvecdef} and \eqref{Gdef}.
Then \eqref{opeqdiscrete} can be rewritten as
\begin{equation}
\label{matrixeq}
\ba \vec{\mathbf{u}} = \vec{\mathbf{f}}.
\end{equation}
Multiplying both side of  \eqref{matrixeq} by $\ba^*$ we get the normal equation
\begin{equation}
\label{noreq}
\L\ba^* \ba\R \vec{\mathbf{u}} = \ba^* \vec{\mathbf{f}}.
\end{equation}
Note that $\ba^* \ba$ is self--adjoint and  positive--definite by the hypothesis.
Note also that $\Ker(\ba^*)$ is orthogonal to $\Ran(\ba)$ and \eqref{Gdecomp} implies that $\Ker(\ba^*)=\Ker(\ba)$.
This and the invertibility of $\ba$ on its range implies that $\ba^* \ba$  is
boundedly invertible on $\Ran({\ba})$.
Therefore, the solution of \eqref{matrixeq} is equivalent to the solution of \eqref{noreq}.

Note that, for $0 < \alpha^* < 2/\lambda_{\max}$, the operator
\begin{equation}
\label{Bdef}
\mathbf{B}:= \alpha^* \sum_{n=0}^\infty (\id -\alpha^* \ba^* \ba)_{|\Ran(\ba)}^n.
\end{equation}
is well--defined and bounded on $\Ran(\ba)$, since $\rho(\alpha^*) := 
\|\left(\id - \alpha^* \ba^* \ba\right)_{|\Ran(\ba)}\|_2 = 
\max\{\alpha^* \lambda_{\max}-1,1- \alpha^* \lambda_{\min}\}< 1$, where $\lambda_{\min}:=\|(\ba^* \ba_{| \Ran(\ba)})^{-1}\|_2$.
The function $\rho$ is minimal at $\alpha^*:=2/(\lambda_{\max}+ \lambda_{\min})$.
Moreover,
\begin{equation}
\label{BGidentity}
\mathbf{B}\circ(\ba^* \ba_{|\Ran(\ba)}) = (\ba^* \ba) \circ \mathbf{B}_{|\Ran(\ba)} =\id_{|\Ran(\ba)}.
\end{equation}
Since $\ba (\id -  \mathbf{P} )=0$,
\begin{equation}
\label{uvecprop}
\ba \vec{\mathbf{u}} =  \ba \mathbf{P} \vec{\mathbf{u}} = \vec{\mathbf{f}}.
\end{equation}
Therefore $ \mathbf{P} \vec{\mathbf{u}} \in \Ran(\ba)$ is the unique solution of \eqref{matrixeq} in $\Ran(\ba)$ and by \eqref{BGidentity}
\begin{equation}
\label{solution}
 \mathbf{P} \vec{\mathbf{u}} =  \mathbf{B} \ba^* \vec{\mathbf{f}}.
\end{equation}
By construction
\begin{eqnarray*}
  \langle f,f_m\rangle_{\bV'\times\bV}&=&\langle\tilde F^*Ff,f_m\rangle_{\bV'\times\bV}\\
  &=&\langle\tilde F^*D_\bV^*\vec{\mathbf f},f_m\rangle_{\bV'\times\bV}\\
  &=&\langle\tilde F^*D_\bV^*\ba\mathbf P\vec{\mathbf u},f_m\rangle_{\bV'\times\bV}\\
  &=&\langle A F^*D_\bV^{-1}\mathbf P\vec{\mathbf u},f_m\rangle_{\bV'\times\bV},\quad m\in\cN,
\end{eqnarray*}
so that $u= F^* D_{\bV}^{-1} \mathbf{P} \vec{\mathbf{u}}$ solves \eqref{basicproblem}.
\end{proof}

\subsubsection{Numerical realization} \label{section:numerrealiz}
Now we turn to the numerical treatment of \eqref{matrixeq}. Due to Theorem \ref{uvectheorem}, the computation
of $\vec \uu$ solving \eqref{matrixeq} amounts to an application of the following damped Richardson iteration
\begin{equation}
\label{richardson}
\vec{\mathbf{u}}^{(i+1)}=\vec{\mathbf{u}}^{(i)}-\alpha^* \ba^* (\ba\vec{\mathbf{u}}^{(i)}-\vec{\mathbf f}),\quad
i \in \NN_0,
\quad {\vec \uu}^{(0)}={\bf 0}.
\end{equation}
Certainly  this iteration cannot be practically realized for
infinite vectors ${\vec \uu}^{(i)}$, $i \in \NN_0$. To avoid this
problem, we make use of the following procedures  (see
\cite{CDD1,CDD2,CDD3,S}  for details on their analysis and
numerical realization) :
\begin{itemize}
\item $\mathbf{RHS}[\varepsilon, \vec{\mathbf{f}}] \rightarrow \vec{\mathbf{f}}_
\varepsilon$: determines for
$\vec{\mathbf{f}} \in \ell_2(\cN)$ a vector $\vec{\mathbf{f}}_\varepsilon 
\in\ell_0(\cN)$ such that
\begin{equation}
\label{RHSprop}
\|\vec{\mathbf{f}} - \vec{ \mathbf{f}}_\varepsilon \|_{\ell_2(\cN)} \leq \varepsilon;
\end{equation}
\item $\mathbf{APPLY}[\varepsilon,\mathbf{A},\vec{ \mathbf{v}}] \rightarrow 
\vec{ \mathbf{w}}_\varepsilon$: determines for a bounded linear operator
$\ba$ on  $\ell_2(\cN)$ and for $\vec{ \mathbf{v}} \in\ell_0(\cN)$ 
a vector $\vec{ \mathbf{w}}_\varepsilon \in \ell_0(\cN)$ such that
\begin{equation}
\label{APPLYprop}
\| \mathbf{A} \vec{ \mathbf{v}} - \vec{ \mathbf{w}}_\varepsilon \|_{\ell_2(\cN)} \leq \varepsilon;
\end{equation}
\item $\mathbf{COARSE}[\varepsilon, \vec{ \mathbf{v}}] \rightarrow \vec{\mathbf{v}}_
\varepsilon$: determines for
$\vec{ \mathbf{v}} \in\ell_0(\cN)$ a vector $\vec{ \mathbf{v}}_\varepsilon \in\ell_0(\cN)$ such that
\begin{equation}
\label{COARSEprop}
\|\vec{ \mathbf{v}} - \vec{ \mathbf{v}}_\varepsilon \|_{\ell_2(\cN)} \leq \varepsilon.
\end{equation}
\end{itemize}

We  discuss in  more details further properties of the routines
$\mathbf{RHS}$, $\mathbf{APPLY}$ and $\mathbf{COARSE}$ in Section
6, where we study the complexity and the computational cost
required to approximate the solution of the original problem up to
some prescribed tolerance.
\\

Let $\rho:=\rho(\alpha^*)$ be as in the proof of Theorem \ref{uvectheorem}. We can define the following inexact version of the damped Richardson
iteration \eqref{richardson}:

\begin{alg}
\label{alg1}
\end{alg}
\fbox{
\begin{minipage}{10.5cm}
 \noindent $\mathbf{SOLVE}[\epsilon, \ba, \vec{\mathbf{f}}]
\rightarrow
\vec{\mathbf{v}}_\epsilon$:\\
\noindent Let $\theta < 1/3$ and $K \in \mathbb{N}$ be fixed such
that $3
\rho^K < \theta$.\\
$j:=0,\;\vec{\mathbf{v}}^{(0)}:=0,\;\epsilon_0:=\|\ba_{|\Ran(\ba)}^{-1}\|\|
\vec{\mathbf{f}} \|_{\ell_2(\cN)}$\\
While $\epsilon_j > \epsilon$ do\\
\indent \hspace{0.5cm} $j:= j+1$\\
\indent \hspace{0.5cm} $\epsilon_j := 3 \rho^K \epsilon_{j-1}/\theta$\\
\indent \hspace{0.5cm} ${\vec{\mathbf{g}}}^{(j)}:=
\mathbf{RHS}[\frac{\theta \epsilon_j}{12 \alpha^* K \| \mathbf
A^*\|},
\vec{\mathbf{f}}]$\\
\indent \hspace{0.5cm} ${\vec{\mathbf{f}}}^{(j)}:=
\mathbf{APPLY}[\frac{\theta \epsilon_j}{12 \alpha^* K}, \ba^*,
\vec{\mathbf{g}}^{(j)}]$\\
\indent \hspace{0.5cm} $\vec{\mathbf{v}}^{(j,0)}:=
\vec{\mathbf{v}}^{(j-1)}$\\
\indent \hspace{0.5cm} For $k=1,...,K$ do \\
\indent \hspace{0.5cm} \indent \hspace{0.5cm}
$\vec{\mathbf{w}}^{(j,k-1)}:=\mathbf{APPLY}[\frac{\theta
\epsilon_j}{12
\alpha^* K \| \mathbf A^*\|}, \ba, \vec{\mathbf{v}}^{(j,k-1)}]$\\
\indent \hspace{0.5cm}  \indent \hspace{0.5cm}
$\vec{\mathbf{v}}^{(j,k)}:= \vec{\mathbf{v}}^{(j,k-1)} - \alpha^*
\left (\mathbf{APPLY}[\frac{\theta \epsilon_j}{12 \alpha^* K},
\ba^*, \vec{\mathbf{w}}^{(j,k-1)}] - \vec{\mathbf{f}}^{(j)} \right)$\\
\indent \hspace{0.5cm} od \\
\indent \hspace{0.5cm}
$\vec{\mathbf{v}}^{(j)}:=\mathbf{COARSE}[(1-\theta)\epsilon_j,\vec{\mathbf{v}}^{(j,K)}]$\\
od \\
$\vec{\mathbf{v}}_\epsilon := \vec{\mathbf{v}}^{(j)}$.
\end{minipage}
}
\\

The parameter $\theta$ plays an important role in  complexity 
estimates (given in Section 6) for $\mathbf{COARSE}$.\\

The proof of the convergence of Algorithm \ref{alg1} is analogous to that of 
\cite[Proposition 2.1]{S} and of \cite[Theorem 4.2]{DFR} except for the fact that 
here we make use of the damped Richardson iteration  in \eqref{richardson} on normal 
equations, due to the non-symmetry of $a$. Nevertheless, the following result holds. 

\begin{tm}
\label{convergence}
Under assumptions  of Theorem \ref{uvectheorem}, let $\vec{\mathbf{u}}\in\ell_2(\cN)$ be a solution of \eqref{matrixeq}. Then $\mathbf{SOLVE}
[\epsilon,\ba,\vec{\mathbf{f}}]$ produces finitely supported vectors $\vec{\mathbf{v}}^{(j,K)},
\vec{\mathbf{v}}^{(j)}, \vec{\mathbf v}_\epsilon$ such that
\begin{equation}
\label{SOLVEprop}
\bigl\|\mathbf P(\vec{\mathbf u}-\vec {\mathbf{v}}^{(j)})\bigr\|_{\ell_2(\cN)}\le\epsilon_j,\quad j\in \NN_0.
\end{equation}
In particular,  
\begin{equation}
\label{SOLVEconv}
\|u-F^*D_\bV^{-1}\vec{\mathbf v}_\epsilon\|_{\bV}\le\|F^*\|\|D_\bV^{-1}\|\epsilon.
\end{equation}
Moreover, it holds that
\begin{equation}
\label{SOLVEhelp}
\bigl\|\mathbf P\vec{\mathbf u}-(\id-\mathbf P)\vec{\mathbf{v}}^{(j-1)}-\vec{\mathbf{v}}^{(j,K)}\bigr\|_{\ell_2(\cN)}\le\frac{2\theta\epsilon_j}{3},
\quad j\ge 1.
\end{equation}
\end{tm}

Of course, the numerical implementation of the damped Richardson
iteration on the normal equations \eqref{noreq} might exhibit a low
convergence rate if the relaxation parameter $\alpha^*$ is small. To
improve the efficiency of the proposed scheme, the generalizations
of Algorithm \ref{alg1} towards, e.g., (conjugate) gradient
iterations as suggested in  \cite{DVU,W} are now a matter of
investigation, see also \cite{DFRSW}.

\begin{figure}
    \centering
    \includegraphics[width=14cm]{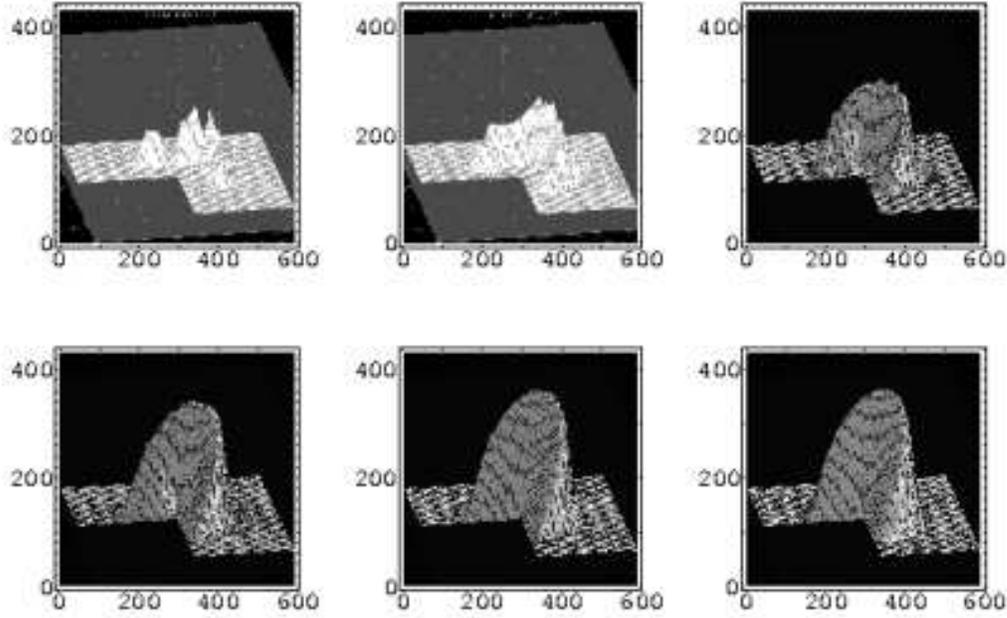}
    \caption{Procedure {\bf SOLVE} applied to the Poisson equation on the L--shaped domain to approximate the solution $u(r,\theta):=\zeta(r)r^{2/3}\sin(\frac{2}{3}\theta)+\tilde u(r,\theta)$, 
    $\zeta\in C^{\infty}(\Omega)$ is a suitable truncation function corresponding to the reentrant corner 
    of the L--shaped domain, and $\tilde u \in H^2(\Omega) \cap H^1_0(\Omega)$. Overlapping rectangular patches and aggregate wavelet frames are used for the discretization. See \cite{DFR,DFRSW,S,W} for more details 
on the numerical implementation. Approximations are shown for the case of
piecewise linear wavelet frame elements.}
    \label{fig:exact solution and rhs 2D}
\end{figure}

\begin{figure}
  \centering
  \includegraphics[width=7cm]{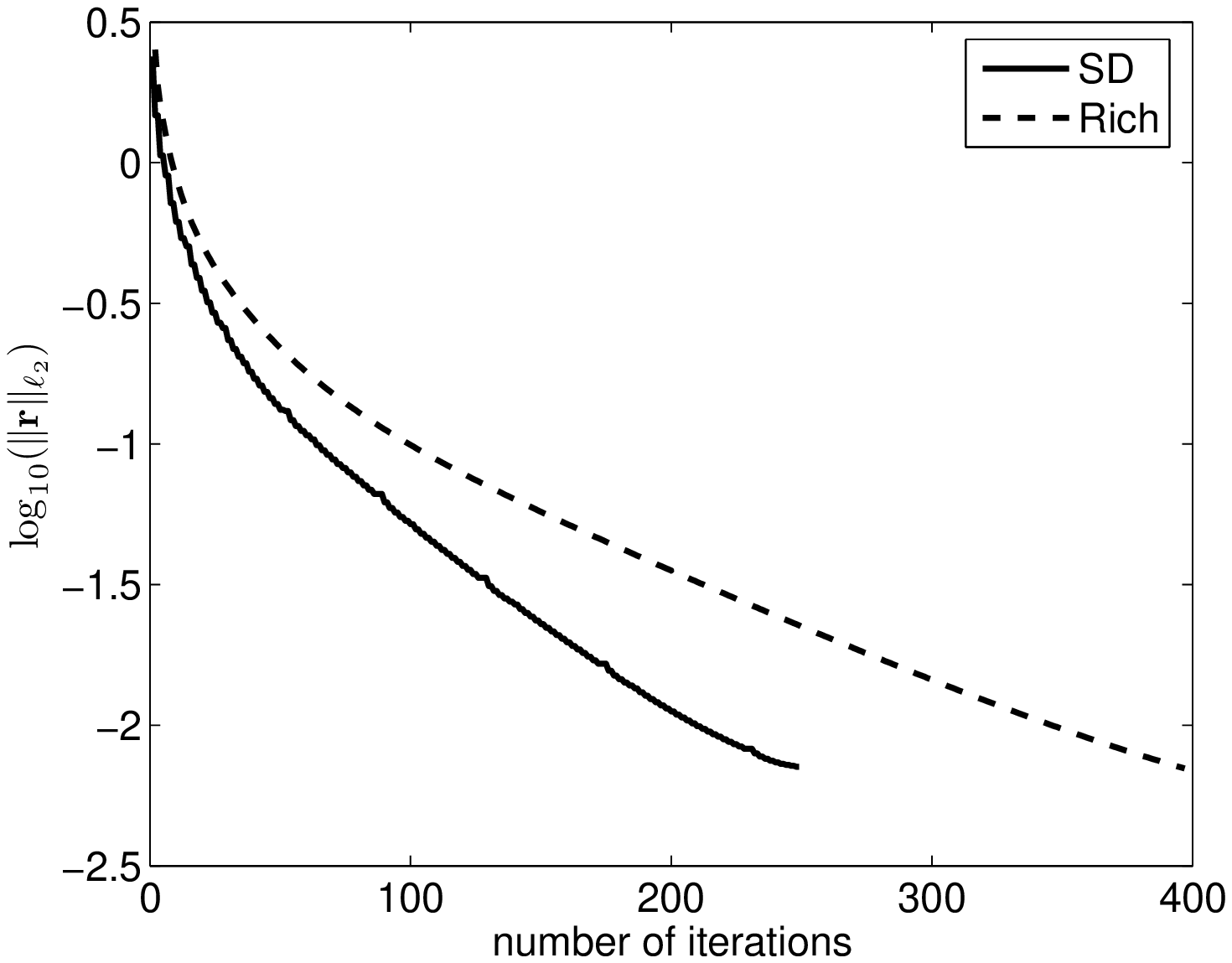}
  \hspace{0.5cm}
  \includegraphics[width=7cm]{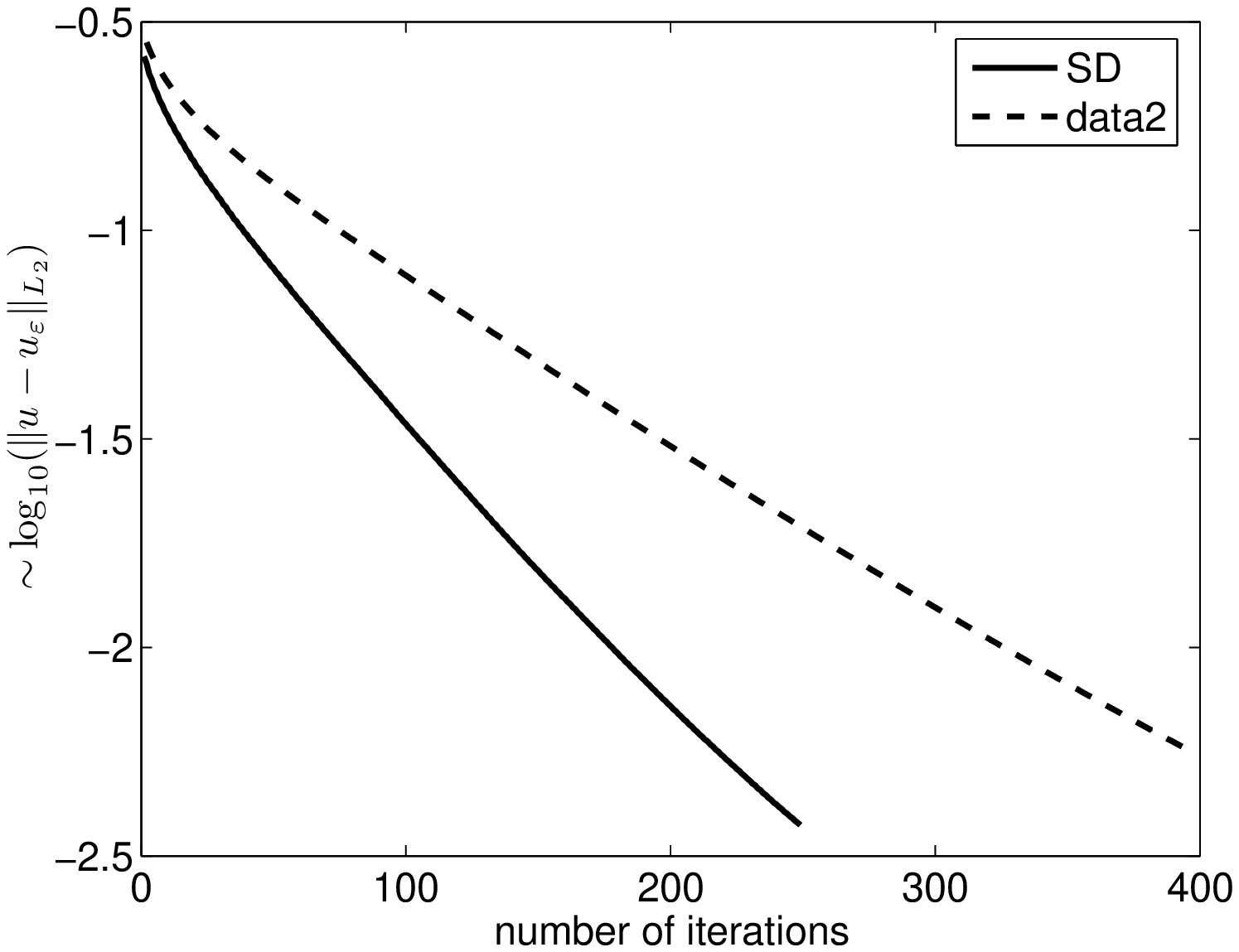}
  \caption{Reduction of the residual $\ell_2$--norm and the $L_2$--error for {\bf SOLVE} (dotted line) and
    for the adaptive steepest descent  (solid line). We refer to \cite{DFRSW} for more detailed numerical tests where optimality is shown and discussed.}
  \label{fig:comp sd rich}
\end{figure}

\section{Numerical realization of the fixed point iteration}

Now we have at hand the major building blocks needed to formulate  an implementable fixed point iteration.
In this section we show how the problem in \eqref{abfor} can be discretized and how Algorithm \ref{alg1} can be used 
to implement the
fixed point iteration in \eqref{fixpt}. 

We denote
\begin{itemize}
\item[i)] $\mathbf{A}:=(D_{\bV}^*)^{-1} F A F^* D_{\bV}^{-1} : \ell_2(\cN) \rightarrow \ell_2 (\cN)$, bounded and boundedly invertible
 on its range $\Ran(\ba) = \Ran((D_{\bV}^*)^{-1} F)$ ;
\item[ii)] $\vec{\bl}:=(D_{\bV}^*)^{-1}F \ell \in \Ran(\ba) \subset \ell_2(\cN)$;
\item[iii)] $\bA1(\cdot):=(D_{\bV}^*)^{-1} F A_1( F^* D_{\bV}^{-1} \cdot, F^*
  D_{\bV}^{-1} \cdot): \ell_2(\cN) \rightarrow \Ran(\ba) \subset \ell_2(\cN)$;
\item[iv)] $\mathbf{H}(\cdot):=(\ba_{|\Ran{\ba}})^{-1} (\vec{\bl} - \bA1 (\cdot)): \ell_2(\cN) \rightarrow \Ran(\ba)$.
\end{itemize}

\begin{rem} Note that $\Ran(\ba)= \Ran((D_{\bV}^*)^{-1} F)$ implies that
the operators in ii)--iv) map $\ell_2(\cN)$ into $\Ran(\ba)$. By
definition of $\PP$
 we also have $\bH(\PP \vec \vv)=\PP \bH( \vec \vv)=\bH(\vec \vv)$ for
any  $\vec \vv \in \ell_2(\cN)$.
\end{rem}

Then, it is easily verified  that the variational  problem in \eqref{abfor} is
equivalent to the following discrete problem.

\vspace{0.3cm}\noindent {\bf Problem 4:} {\it Find $\vec{\mathbf{u}} \in \Ran(\ba) \subset \ell_2(\cN)$ such that
\begin{equation}
\label{abfor2}
\ba \vec{\mathbf{u}} + \bA1(\vec{\mathbf{u}}) = \vec{\bl},
\end{equation}
or, equivalently, such that $\vec{\mathbf{u}}$ is a fixed point  of $\mathbf{H}$ in $\Ran(\ba)$, i.e
\begin{equation}
\label{abfor3}
\vec{\mathbf{u}} = \mathbf{H}(\vec{\mathbf{u}}), \quad \vec \uu \in \Ran(\ba).
\end{equation}}

Define the closed subset of $\Ran(\ba)$ by
$$
 \BB_r :=\{\vec{\mathbf{u}} \in \Ran(\ba):\|\vec{\mathbf{u}}\|_{\ell_2(\cN)} \leq r\} \quad
 \text{for some} \quad r \in \RR_+.
$$
The next theorem is a discrete analogue of Lemma 2  and Corollary
2 in \cite{CMS}.

\begin{tm}
\label{discrprob}
Under the assumptions and notations specified above, the following statements hold true:
\begin{itemize}
\item[a)] If $0<r < \left( 2 \|(\ba_{|\Ran{\ba}})^{-1}\| \|F\|^3 \| A_1\| \right) ^{-1}$, then
$\mathbf{H}_{|\BB_r}$ is a contraction with Lipschitz constant
$L:= r \left (2 \|(\ba_{|\Ran{\ba}})^{-1}\| \|F\|^3 \| A_1\| \right) <1$;

\item[b)] if $0<r < \left( \|(\ba_{|\Ran{\ba}})^{-1}\| \|F\|^3 \| A_1\| \right) ^{-1}$ and $\| \vec{\bl} \| \leq r \L \|(\ba_{|\Ran{\ba}})^{-1}\|^{-1} - \| A_1\|  \|F\|^3 r\R$ then $\mathbf{H}(\BB_r) \subseteq  \BB_r$ ;
\item[c)] if $\| \vec{\bl} \| < \L 4 \|(\ba_{|\Ran{\ba}})^{-1}\|^2 \|F\|^3 \| A_1\|\R^{-1}$ then
\eqref{abfor3} has a unique solution $\vec{\mathbf{u}} \in \BB_{r^*}$, for some suitable $r^*$ such that
$0<r^*<  \left( 2 \|(\ba_{|\Ran{\ba}})^{-1}\| \|F\| \| A_1\| \right) ^{-1}$.
\end{itemize}
\end{tm}

\begin{proof}
Recall that $D_{\bV}$ is assumed unitary and $\|D_{\bV}\|=\|D_{\bV}^{-1}\|=\|D_{\bV}^*\| \equiv 1$. 
Moreover, since $F$ maps $\bV'$ into $\bV_d'$ both being Hilbert spaces, thus   $\|F\|=\|F^*\|$. Then, for $\vec \uu,\vec \vv \in \BB_r$,
\begin{eqnarray*}
\|\mathbf{H}\vec \uu - \mathbf{H} \vec \vv\| &\leq & 
\|(\ba_{|\Ran{\ba}})^{-1}\| \|\bA1 \vec \uu - \bA1 \vec \vv\|_{\ell_2(\cN)}\\
&=& \|(\ba_{|\Ran{\ba}})^{-1}\| \| (D_{\bV}^*)^{-1} F  \L A_1( F^*  D_{\bV}^{-1} \vec 
\uu, F^*  D_{\bV}^{-1}\vec \uu) -  A_1( F^*  D_{\bV}^{-1} \vec \vv, F^*  D_{\bV}^{-1} \vec \vv)\R\|_{\ell_2(\cN)}\\
&\leq & \|(\ba_{|\Ran{\ba}})^{-1}\|  \|F\| \|  A_1( F^*  D_{\bV}^{-1} \vec \uu, F^*  D_{\bV}^{-1} 
\vec \uu) -  A_1( F^*  D_{\bV}^{-1}\vec \vv, F^*  D_{\bV}^{-1} \vec \vv)\|_{\bV'}\\
&=& \|(\ba_{|\Ran{\ba}})^{-1}\|  \|F\| \|  A_1( F^* D_{\bV}^{-1} (\vec \uu-\vec \vv), F^* 
D_{\bV}^{-1} \vec \uu) -  A_1( F^*  D_{\bV}^{-1} \vec \vv, F^*  D_{\bV}^{-1}(\vec \vv-\vec \uu))\|_{\bV'}\\
&\leq& \|(\ba_{|\Ran{\ba}})^{-1}\|  \|F\| \| A_1\| \|F^*\|^2\|\vec \uu-\vec \vv\|_{\ell_2(\cN)} 
\left(\|\vec \uu\|_{\ell_2(\cN)}+
\|\vec \vv\|_{\ell_2(\cN)} \right) \\
&\leq&  2r  \|(\ba_{|\Ran{\ba}})^{-1}\|  \|F\|^3 \| A_1\|   \|\vec \uu-\vec \vv\|_{\ell_2(\cN)}.
\end{eqnarray*}
Thus, as  $L <1$ we have that $\mathbf{H}$ is a contraction.
To show b), just observe that by definition of $\mathbf{H}$ and the estimate above
\begin{eqnarray*}
\|\mathbf{H}\vec \uu\|_{\ell_2(\cN)} &\leq & \| \ba_{|\Ran{\ba}} ^{-1}\| \L 
\| \vec{\bl}\| + \| A_1\| \|F\|^3 r^2 \R \leq r.
\end{eqnarray*}
c) We have to show that if $\| \vec{\bl} \| < \L 4 \|(\ba_{|\Ran{\ba}})^{-1}\|^2 \|F\|^3 \| A_1\|\R^{-1}$
then there exists $r^*$ with $0<r^*<\left( 2 \|(\ba_{|\Ran{\ba}})^{-1}\| \|F\|^3 \| A_1\| \right) ^{-1}$
such that $\|\vec{\bl}\|< r^* \L \|(\ba_{|\Ran{\ba}})^{-1}\|^{-1} - \| A_1\|  \|F\|^3 r^*\R$. Then, using
parts a) and b) of this Lemma we get that $H|_{\BB_{r^*}}$ is a contractive mapping of $\BB_{r^*}$
into itself and, thus, has a unique fixed point. To see that such an $r^*$ exists, consider
$$h(r)=r \L \|(\ba_{|\Ran{\ba}})^{-1}\|^{-1} - \| A_1\|  \|F\|^3 r\R,$$
a quadratic mapping, and
note that $h$ assumes values from $0$ to $\L 4 \|(\ba_{|\Ran{\ba}})^{-1}\|^2 \|F\|^3 \| A_1\|\R^{-1}$ as
$r$ varies from $0$ to $\left( 2 \|(\ba_{|\Ran{\ba}})^{-1}\| \|F\|^3 \| A_1\| \right) ^{-1}$.
\end{proof}

\begin{cor}
\label{discrsolution} If $\| \vec{\bl} \| < (4 \|(\ba_{|\Ran{\ba}})^{-1}\|^2 \|F\|^3 \|A_1\|)^{-1} $,
then the solution $\vec \uu$
in $\BB_{r^*}$ of \eqref{abfor3} is given by the following discrete fixed point iteration
\begin{eqnarray}
\label{fixpt2}
\vec \uu_{n+1} &=& \mathbf{H} \vec \uu_{n}, \quad n \geq 0, \quad \vec \uu_0 = \mathbf{0} \in \Ran(\ba), \\
\vec \uu&=& \lim_{n \rightarrow \infty} \vec \uu_n.
\end{eqnarray}
\end{cor}

The discrete fixed point iteration cannot be implemented  for
infinite vectors. For this reason we propose the following new
approximating adaptive scheme $\mathbf{FIXPT}$, where the
procedure $\mathbf{SOLVE}$ introduced in Subsection
\ref{section:numerrealiz} replaces the exact computation  of $\vec \uu_{n+1}$ in
\eqref{fixpt2}.

\begin{alg} \label{alg2}
\end{alg}
\fbox{
\begin{minipage}{7cm}
\noindent $\mathbf{FIXPT}[\varepsilon, \ba,\bA1,\vec{\bl}]
\rightarrow
\vec{\mathbf{u}}_\varepsilon$:\\

\noindent $i:=0,\;\vec \vv_0:=0,\;0< \varepsilon_0 < r^*; \
\varepsilon_0 \not = L$;\\ While $\L\varepsilon_i>\frac{\varepsilon_0-L}{\varepsilon_0 \L
\varepsilon_0-L\L\frac{L}{\varepsilon_0}\R^i \R} (\varepsilon -
L^i r^*)\R$ do\\ \indent \hspace{0.5cm} $\varepsilon_{i+1} :=
\varepsilon_0^{i+1}$\\ \indent \hspace{0.5cm} $\vec \vv_{i+1}
=\mathbf{SOLVE}[\varepsilon_{i+1},\ba,\vec{\bl}\mathbf{
-\bA1(}\vec \vv_{i}\mathbf{)]}$\\ \indent \hspace{0.5cm} $i:=
i+1$\\ \noindent od\\ $\vec{\mathbf{u}}_\varepsilon := \vec
\vv_i$. 
\end{minipage}
}
\\

\begin{rem}
The Algorithm \ref{alg2} is a perturbation of the exact fixed point iteration creating  sequences
 $\{\vec \vv_i\}_{i\in \NN_0}$ of finite vectors. In general such vectors will not belong to $\Ran(\ba)$ and there
 is not much hope that $\lim_{n \rightarrow \infty} \vec \vv_n = \vec \uu$. However, one might try to see whether
 $\lim_{n \rightarrow \infty} \mathbf{P} \vec \vv_n = \vec \uu$. A priori it might even happen that, because of an
 accumulation of the perturbation errors, there exists some $n \in \NN$ large enough such that $\mathbf{P}\vec \vv_n \notin \BB_{r^*}$!
 In order to show the convergence of  Algorithm \ref{alg2} to the solution of  problem \eqref{abfor3} we prove the following Lemma.
\end{rem}

\begin{lemma} \label{bound}If $\| \vec{\bl} \| <  (4 \|(\ba_{|\Ran{\ba}})^{-1}\|^2 \|F\|^3 \|A_1\|)^{-1}$ 
and if the quantities
\begin{eqnarray*}
\mathcal{E}_n&:=&\sum_{k=2}^{n+1} \varepsilon_k + \L 1+L \R \sum_{h=0}^{n-3} \sum_{k=3}^{n-h} 
\varepsilon_k L^{n-h-k}+\L \varepsilon_1 +
 L(\varepsilon_0+ \|(\ba_{|\Ran(\ba)})^{-1}\| \|\vec{\bl}\|) \R  \sum_{k=0}^{n-1} L^{k}  \\
&+&  \varepsilon_0+\|\ba_{|\Ran(\ba)}^{-1}\| \|\vec{\bl}\| \leq r^*, \text{ for all } n \in \NN_0,
\end{eqnarray*}
with $L$ and $r^*$ as in Theorem \ref{discrprob} a) and
c), respectively, then the sequence $\{\mathbf{P} \vec \vv_i\}_{i
\in \mathbb{N}}$ resulting from the application of Algorithm
\ref{alg2} all lies in $\BB_{r^*}$.
\end{lemma}

\begin{proof}
We show by induction on $n$ that
\begin{equation}
\label{iter}
\| \mathbf{P}(\vec \vv_{n+1}-\vec \vv_{n})\|_{\ell_2(\cN)} \leq \varepsilon_{n+1} + \L 1+L \R  
\sum_{k=3}^n \varepsilon_k L^{n-k}
   +L^{n-1}(\varepsilon_1 + L\|\mathbf{P} \vec \vv_1\|_{\ell_2(\cN)}).
\end{equation}
and
\begin{eqnarray}
\label{iter1}
\| \mathbf{P}(\vec \vv_{n+1})\|_{\ell_2(\cN)} &:=&\sum_{k=2}^{n+1} \varepsilon_k +
\L 1+L \R \sum_{h=0}^{n-3} \sum_{k=3}^{n-h} \varepsilon_k L^{n-h-k}  \\  \nonumber
 &+&\L \varepsilon_1 + L(\varepsilon_0+ \|(\ba_{|\Ran(\ba)})^{-1}\| \|\vec{\bl}\|) \R  
 \sum_{k=0}^{n-1} L^{k} 
 +  \varepsilon_0+\| \ba_{|\Ran(\ba)} ^{-1}\| \|\vec{\bl}\|
\end{eqnarray}
Set $n=1$. We get 
\begin{eqnarray*}
 \|\mathbf{P} \vec \vv_1\|_{\ell_2(\cN)} = 
 \| \mathbf{P} \L \mathbf{SOLVE[\varepsilon_0,\ba,} \vec{\bl}\mathbf{]}-\mathbf{H}(0)\R
 + \mathbf{H}(0)\|_{\ell_2(\cN)} &\leq& \varepsilon_0 + \| \ba_{|\Ran(\ba)} ^{-1} \vec{\bl}\| \\ &\leq& \varepsilon_0 +
 \|(\ba_{|\Ran(\ba)})^{-1} \|\| \vec{\bl}\|\leq r^*.
\end{eqnarray*}
Moreover, since $\mathbf{H}(\vec \vv)=\mathbf{H}(\mathbf{P} \vec
\vv)= \PP \bH(\vec \vv)$, we get
$$\mathbf{P} \L
\mathbf{SOLVE[\varepsilon_2,\ba,}\vec{\mathbf{l}}\mathbf{-\bA1(}\vec
\vv_1\mathbf{)]} \R- \mathbf{H}(\vec \vv_1)= \mathbf{P} \L
\mathbf{SOLVE[\varepsilon_2,\ba,}\vec{\mathbf{\bl}}\mathbf{-\bA1(}\vec{\mathbf{\vv_1}}\mathbf{)]} \R -
\mathbf{P} \;\mathbf{H}(\vec \vv_1).$$
 By \eqref{SOLVEprop} it holds that
$\|\mathbf{P} \ \mathbf{SOLVE}[\varepsilon_{n+1},\ba,\vec{\bl}-\bA1(\vec
\vv_n)] - \mathbf{P} \ \mathbf{H}(\vec \vv_n)\| \leq\varepsilon_{n+1}$ and,
using the assumption on $\|\vec{\bl}\|$ ensuring that $\mathbf H$ is a contraction with the 
Lipschitz constant $L$, we get
\begin{eqnarray*}
\| \mathbf{P}(\vec \vv_{2}-\vec \vv_{1})\|_{\ell_2(\cN)} &=& \left\|\mathbf{P} \L 
\mathbf{SOLVE}[\varepsilon_2,\ba,\vec{\bl}-\bA1(\vec \vv_1)] \R-
\mathbf{P} \L\mathbf{SOLVE}[\varepsilon_1,\ba,\vec{\bl}-\bA1(\vec \vv_0)] \R \right\|_{\ell_2(\cN)} \\
&=& \| \mathbf{P} \L \mathbf{SOLVE}[\varepsilon_2,\ba,\vec{\bl}-\bA1(\vec \vv_1)]- 
\mathbf{H}(\vec \vv_1) \R+
 \mathbf{H}(\vec \vv_1)-\mathbf{H}(\vec \vv_0)  \\
&+& \mathbf{H}(\vec \vv_0)- \mathbf{P} \L\mathbf{SOLVE}
[\varepsilon_1,\ba,\vec{\bl}-\bA1(\vec \vv_0)]\R \|_{\ell_2(\cN)} \\
&\leq& \varepsilon_2 +\varepsilon_1+L\|\mathbf{P} \vec \vv_1\|_{\ell_2(\cN)} .
\end{eqnarray*}
We assume now that all $\mathbf{P} \vec \vv_1,...,\mathbf{P}\vec \vv_n \in \BB_{r^*}$ and that formulas
\eqref{iter}  and \eqref{iter1} are valid  replacing $n$ with $n-1$.
Then
\begin{eqnarray*}
\| \mathbf{P}(\vec \vv_{n+1}-\vec \vv_{n})\|_{\ell_2(\cN)}  &\leq&  
\varepsilon_{n+1} +\varepsilon_n+L\|\mathbf{P} (\vec \vv_n -\vec \vv_{n-1})\|_{\ell_2(\cN)} \\
&\leq& \varepsilon_{n+1} +\varepsilon_n+L \L \varepsilon_{n} + \L 1+L \R  \sum_{k=3}^{n-1} 
\varepsilon_k L^{n-1-k}   +L^{n-2}(\varepsilon_1 + L\|\mathbf{P} \vec \vv_1\|_{\ell_2(\cN)} ) \R\\
&=& \varepsilon_{n+1} + \L 1+L \R  \sum_{k=3}^n \varepsilon_k L^{n-k}   +L^{n-1}(\varepsilon_1 + 
L\|\mathbf{P} \vec \vv_1\|_{\ell_2(\cN)} ).
\end{eqnarray*}
Using the above estimate for $\| \mathbf{P}(\vec \vv_{n+1}-\vec \vv_{n})\|_{\ell_2(\cN)} $, the triangular inequality,
i.e.
$$\| \mathbf{P} \vec \vv_{n+1} \|_{\ell_2(\cN)}   \le \| \mathbf{P}(\vec \vv_{n+1}-\vec \vv_{n})\|_{\ell_2(\cN)} 
+\| \mathbf{P} \vec \vv_{n} \|_{\ell_2(\cN)} $$
and \eqref{iter1} for $\mathbf{P} \vec \vv_{n } $ we have
immediately \eqref{iter1} for $\mathbf{P} \vec \vv_{n+1} $.
Thus, by hypothesis on $\mathcal{E}_n$ we obtain $ \|\mathbf P
\vec \vv_{n+1}\|_{\ell_2(\cN)}  \leq r^*$. This implies by induction that $
\mathbf{P} \vec \vv_{n} \in \BB_{r^*}$ for all $n \in \mathbb{N}$.
\end{proof}

\begin{rem} Note that the assumption on $\mathcal{E}_n$ in the
above Theorem is not restrictive as soon as the data are small. In fact, it holds that
\begin{itemize}

\item $\displaystyle{\sum_{k=2}^{n+1} \varepsilon_k=\sum_{k=0}^{n+1} \varepsilon_0^k-1-\varepsilon_0
 =\frac{1-\varepsilon_0^{n+2}}{1-\varepsilon_0}-1-\varepsilon_0 \leq \frac{\varepsilon_0^2}{1- \varepsilon_0}
 \rightarrow 0}$ for $\varepsilon_0 \rightarrow 0$;
\item the map
$$n \mapsto  \sum_{h=0}^{n-3} \sum_{k=3}^{n-h} \varepsilon_k L^{n-h-k}  \leq C(\varepsilon_0,L),$$
where $C(\varepsilon_0,L) \rightarrow 0$ for $\varepsilon_0 \rightarrow 0$, uniformly with respect to $n$;

\item $ \displaystyle{\L\varepsilon_1+L(\varepsilon_0+ \| \ba_{|\Ran(\ba)}^{-1}\| \|\vec{\bl}\| 
\R  \sum_{k=0}^{ \infty} L^{k}  +
 (\varepsilon_0+\| \ba_{|\Ran(\ba)}^{-1}\| \|\vec{\bl}\|)}$ is small as one wants whenever $\varepsilon_0$ and $\|\vec{\bl}\|$
 are small enough. Note
that $\|\vec \bl\| \lesssim \|F\| \|\ell\|$, where $\|\ell\|$, at
least for the MHD problem discussed above, depends on the size of
the norms of the forcing terms and boundary data in
\eqref{momentum.eq}--\eqref{elec potential} by means of formula \eqref{bound:ell}.
Unfortunately for nonlinear problems in fluid dynamics existence and uniqueness of the solutions are ensured only under smallness assumptions as in \eqref{bound:ell}, so that we cannot expect that numerical implementations can be less restrictive. Nevertheless, such restrictions represent as usual the worst case analysis.
\end{itemize}
\end{rem}

We are now ready to show our main convergence result.
\begin{tm}
\label{converg}
 If $\| \vec{\bl} \| <   (4 \|(\ba_{|\Ran{\ba}})^{-1}\|^2 \|F\|^3 \|A_1\|)^{-1}$ 
and $\mathcal{E}_n \leq r^*$ for all $n \in \NN$,
 then
\begin{equation}
\mathbf{H}(\vec \uu)=\vec \uu = \lim_{n \rightarrow \infty} \mathbf{P} \vec \vv_n =\mathbf{P}\L\mathbf{FIXPT}[0, \ba,\bA1,\vec{\bl}]\R,
\end{equation}
where $\{\vec \vv_i\}_{i \in \mathbb{N}}$ are obtained by Algorithm \ref{alg2}.
After $n$ iterations of Algorithm \ref{alg2} one gets
\begin{equation}
\label{iteration1}
\|\mathbf{P} \vec \vv_{n+1} - \vec \uu\|_{\ell_2(\cN)}  
\leq \varepsilon_0 \frac{\varepsilon_0^n (\varepsilon_0- L \L\frac{L}{\varepsilon_0}
\R^n )}{\varepsilon_0-L}+  L^n\|\vec \uu\| \leq \varepsilon_0 \frac{\L\varepsilon_0^n 
(\varepsilon_0- L \L\frac{L}{\varepsilon_0}
\R^n \R}{\varepsilon_0-L}+  L^n r^*.
\end{equation}
Therefore, for $\varepsilon>0$ one has
\begin{equation}
\label{estimation1}
\| \vec \uu-\mathbf{P} \L\mathbf{FIXPT}[\varepsilon, \ba,\bA1,\vec{\bl}]\R\| \leq \varepsilon,
\end{equation}
and the number $N$ of iterations to achieve the accuracy
$\varepsilon>0$ is estimated by
\begin{equation}
\label{numiter} N \sim - \log(\varepsilon).
\end{equation}
\end{tm}

\begin{proof}
First we want to show that $\{\mathbf{P} \vec \vv_i\}_{i \in
\mathbb{N}}$ is a Cauchy sequence in $\BB_{r^*}$. To do that, we
use the estimation \eqref{iter} and get
\begin{eqnarray*}
\|\mathbf P (\vec \vv_{n+m}- \vec \vv_{n})\|_{\ell_2(\cN)} &\leq&
\sum_{k=n+1}^{n+m} \varepsilon_k + \L 1+L \R   \sum_{h=1}^{m} \sum_{k=3}^{n+m-h} \varepsilon_k 
L^{n+m-h-k} \\
&+&\L \varepsilon_1 + L(\varepsilon_0+ \|(\ba_{|\Ran(\ba)})^{-1}\| \|\vec{\bl}\|) \R 
 \sum_{k=n-1}^{n+m-2} L^{k} .
\end{eqnarray*}
A straight forward computation yields that the expression on the
right in the estimate above goes to zero as $n$ goes to infinity.
Therefore, for any $\varepsilon>0$ there exists $n \in \NN$ large
enough such that for all $m \in \NN$, $m
>n$ it holds that $\|\mathbf P (\vec \vv_{n+m}- \vec \vv_{n})\|_{\ell_2(\cN)}  \leq
\varepsilon$. Due to Lemma  \ref{bound} the sequence $\{\mathbf{P}
\vec \vv_i\}_{i \in \mathbb{N}} \in \BB_{r^*}$. Therefore,
$\{\mathbf{P} \vec \vv_i\}_{i \in \mathbb{N}}$ is a Cauchy
sequence in $\BB_{r^*}$ and then (because $\BB_{r^*}$ is a closed
subset) there exists a unique $\vec \vv \in \BB_{r^*}$ such that
$\displaystyle{\vec \vv=  \lim_{n \rightarrow \infty} \mathbf{P}
\vec \vv_n}$. Moreover, we have
\begin{eqnarray*}
\mathbf{P} \vec \vv_{n+1} &=& \mathbf{P}\L \mathbf{SOLVE}[\varepsilon_{n+1},\ba,\vec{\bl}-\bA1(\vec \vv_n)]\R \\
&=& \L \mathbf{P}
\L\mathbf{SOLVE}[\varepsilon_{n+1},\ba,\vec{\bl}-\bA1( \vec
\vv_n)]\R - \mathbf{H}(\mathbf{P}\vec \vv_n)\R +
\L\mathbf{H}(\mathbf{P}\vec \vv_n)-\mathbf{H}(\vec \uu) \R +\vec
\uu .
\end{eqnarray*}

This implies that
\begin{eqnarray*}
\|\mathbf{P} \vec \vv_{n+1} - \vec \uu\|   &=& 
\| \mathbf{P}\L \mathbf{SOLVE}[\varepsilon_{n+1},\ba,\vec{\bl}-\bA1(\vec \vv_n)]\R -
\mathbf{H}(\mathbf{P}\vec \vv_n) \|_{\ell_2(\cN)} + 
\| \mathbf{H}(\mathbf{P}\vec \vv_n)-\mathbf{H}(\vec \uu)\|_{\ell_2(\cN)}  \\
&\leq& \varepsilon_{n+1} +L\|\mathbf{P}\vec \vv_n-\vec \uu\|_{\ell_2(\cN)}  
\leq \varepsilon_{n+1} +L( \varepsilon_n+L
\|\mathbf{P}\vec \vv_{n-1}-\vec \uu\|_{\ell_2(\cN)} )\\
&\leq& \varepsilon_0 \sum_{k=0}^n \varepsilon_0^{n-k} L^k+ L^n\|\vec \uu\|_{\ell_2(\cN)} \\
&=& \varepsilon_0 \frac{(\varepsilon_0^{n+1}- L^{n+1})}{\varepsilon_0-L}+  L^n\|\vec \uu\|_{\ell_2(\cN)}  \\
&=& \varepsilon_0 \frac{ \varepsilon_0^n \L\varepsilon_0- L
\L\frac{L}{\varepsilon_0}\R^n \R}{\varepsilon_0-L}+  L^n\|\vec
\uu\|_{\ell_2(\cN)}  \rightarrow 0, \quad n \rightarrow \infty.
\end{eqnarray*}
Therefore,  one has $\mathbf{H}(\vec \uu)=\vec \uu = \lim_{n \rightarrow \infty} \mathbf{P} \vec \vv_n =\mathbf{P}\L\mathbf{FIXPT}
[0, \ba,\bA1,\vec{\bl}]\R$.
Moreover, since $\|\vec \uu\|_{\ell_2(\cN)}  \leq r^*$, then the second inequality of \eqref{iteration1} is valid and one has
\begin{equation}
\label{estimation2}
\varepsilon_0 \frac{ \varepsilon_0^n (\varepsilon_0- L \L\frac{L}{\varepsilon_0}\R^n )}{\varepsilon_0-L}+  L^n r^* \leq \varepsilon
\end{equation}
if and  only if
\begin{equation}
\label{estimation3}
\varepsilon_n \leq \frac{\varepsilon_0-L}{\varepsilon_0 \L \varepsilon_0-L\L\frac{L}{\varepsilon_0}\R^n \R} (\varepsilon - L^n r^*),
\end{equation}
which is the criterion to exit the loop in Algorithm \ref{alg2}.
This implies immediately \eqref{estimation1}. From \eqref{estimation2} one shows that the necessary number of iterations of achieve accuracy $\varepsilon>0$ is given by \eqref{numiter}.
\end{proof}

Similarly to the proof of formula \eqref{SOLVEconv} (see \cite[Theorem 4.2]{DFR}) one finally
shows the following.

\begin{cor}
\label{converg2}
 If $\| \vec{\bl} \| <   (4 \|(\ba_{|\Ran{\ba}})^{-1}\|^2 \|F\|^3 \|A_1\|)^{-1}$ 
and $\mathcal{E}_n \leq r^*$ for all $n \in \NN$, then
\begin{equation}
u = \sum_{n\in\mathcal{N}} \L D_{\bV}^{-1} \mathbf{FIXPT}[0,
\ba,\bA1,\vec{\bl}]\R_n f_n,
\end{equation}
is a solution in $\bV$ of the fixed point problem $H(u)=u$ and one has that for $\varepsilon>0$
\begin{equation}
\| u -  \sum_{n\in\mathcal{N}} \L D_{\bV}^{-1} \mathbf{FIXPT}[\varepsilon, \ba,\bA1,\vec{\bl}]\R_n f_n \|_{\bV} \lesssim \varepsilon.
\end{equation}
\end{cor}
\begin{proof}
Since $\Ker(F^*D_\bV^{-1})=\Ker(\ba)=\Ker(\mathbf P)$ and $F^*
\tilde{F}=\id$, we finally verify
\begin{eqnarray*}
\|u-F^*D_\bV^{-1}\vec{\mathbf u}_{\varepsilon}\|_{\bV}
&=&\bigl\|F^*(\tilde F u-D_\bV^{-1}\vec{\mathbf u}_{\varepsilon})\|_{\bV}\\
&=&\bigl\|F^*D_\bV^{-1} (\mathbf P \vec{\mathbf u}-\vec{\mathbf u}_{\varepsilon})\bigr\|_{\bV}\\
&=&\bigl\|F^*D_\bV^{-1} \mathbf P (\vec{\mathbf u}-\vec{\mathbf u}_{\varepsilon})\bigr\|_{\bV}\\
&\le&\|F^*\|\|D_\bV^{-1}\|\bigl\|\mathbf P(\vec{\mathbf u}-\vec{\mathbf u}_{\varepsilon})\bigr\|_{\ell_2(\cN)}.
\end{eqnarray*}
The claim follows from the above Theorem.
\end{proof}

\begin{rem}
The fact that at each iteration we compute the approximation up to a perturbation/tolerance 
$\varepsilon_i$ also means that the scheme is stable.
In other words, not only can $\varepsilon_i$  be interpreted as the numerical approximation accuracy
we achieve at each step, but also as the error tolerance we can afford, without spoiling convergence. Moreover, the 
scheme is {\it fully adaptive} in the sense that the iterations are enforced (by the use of 
suitable implementation of $\mathbf{COARSE}$, see below) to work only with  minimal number of relevant 
quantities (frame coefficients), in order to keep the prescribed accuracy--complexity balance.
\end{rem}

\section{Quasi--optimal complexity of the algorithm}
In this section we present the complexity analysis of Algorithm \ref{alg2}.
From Theorem \ref{converg}, in particular from formula
\eqref{numiter}, we already know that to achieve a prescribed
accuracy $\varepsilon>0$ one needs to execute $N \sim -
\log(\varepsilon)$ iterations. Therefore, having an estimation of
the cost of each iteration, the asymptotic analysis of the
complexity of the suggested adaptive scheme can be done. This
is one of the very interesting theoretical advantages of the
adaptive (wavelet) frame approach, together with the fact that one
can prove both convergence and stability of the adaptive scheme as
shown in the previous section.

Of course, the main ingredient of iterations  in Algorithm
\ref{alg2} is the procedure $\mathbf{SOLVE}$. As
announced in Subsection \ref{section:numerrealiz}, we discuss here
an implementation of such a procedure and study its complexity.
To this end, one should illustrate how the building block
procedures $\mathbf{RHS}$, $\mathbf{APPLY}$, and $\mathbf{COARSE}$
can be implemented and estimate their computational cost.
Therefore, in this section we focus on the main {\it properties
and requirements} of these building blocks so that Algorithms \ref{alg1} and \ref{alg2} have certain complexity.
We refer the interested reader to \cite{CDD2,CDD3,DFR,S} for the
descriptions of procedures $\mathbf{RHS}$, $\mathbf{APPLY}$,
and $\mathbf{COARSE}$ that fulfill the requirements stated below
and  needed for the complexity estimates. The complexity
estimates for even more general algorithms than Algorithm
\ref{alg2} for linear and nonlinear variational problems,
but under the more restrictive assumption that the discretizing
frame $\mathcal{F}$ is a Riesz basis, have been given in
\cite{CDD2,CDD3}. In order to describe the complexity in the more
general case of pure frame discretizations a bit more (technical) effort and
preparation is needed. We start by defining a so--called sparseness
class $\cA^s$ of vectors, $s \in \RR_+$. $\cA^s$ will turn out to be such
that, if $\vec \uu \in \cA^s$, then the size of the support of
$\vec \uu_\varepsilon = \mathbf{FIXPT}[\varepsilon,
\ba,\bA1,\vec{\bl}]$ and the computational cost for obtaining
$\vec \uu_\varepsilon$ can be estimated a priori. 

 An algorithm for computing a finite approximation $\vec \uu_\varepsilon$ of $\vec
\uu$  up to $\varepsilon$, for $\vec \uu$
 given implicitly as a solution of some equation, is called {\it optimal} if $\# \supp(\vec \uu_\varepsilon)$ (the number of elements of the support of $\vec \uu_\varepsilon$) is not asymptotically
 larger (for $\varepsilon \rightarrow 0$) than the same quantity obtained by direct computation of any
 other approximation of $\vec \uu$ using the same tolerance $\varepsilon$, for $\vec
\uu$ being given explicitly. In addition to this,
  the optimality is fully realized if the complexity to compute $\vec \uu_\varepsilon$ does not exceed $\# \supp(\vec \uu_\varepsilon)$ asymptotically
  (for $\varepsilon \rightarrow 0$). In other words, one does want
  the number of algebraic operations to be comparable to the size of what is being computed.
For discretizations by means of Riesz bases optimal algorithms
can be realized, see \cite{CDD1,CDD2,CDD3}, for example. Analogous
algorithms for frames may exhibit ``arbitrarily small reductions'' with
respect to the expected optimality (see the following Theorem 6.1 for the precise statement). 
Since the techniques  for estimating complexity we use are similar to those in \cite[Theorem 3.12]{S}, 
we encounter similar difficulties to achieve the full optimality for $\mathbf{FIXPT}$ when dealing 
with pure frames. We call this situation {\it quasi--optimal}.

An optimal sparseness class is modeled as follows: For given
$s>0$ we define the space 
\begin{equation}
\label{weakelltau} \cA^s_{weak}:=\{\vec{\mathbf
c}\in\ell_2(\cN):\;\|\vec{\mathbf
c}\|_{\cA^s_{weak}}:=\sup_{n\in\mathbb
N}n^{1/2+s}|\gamma_n(\vec{\mathbf c})|<\infty\},
\end{equation}
where $\gamma_n(\vec{\mathbf c})$ is the $n-$th largest
coefficient in modulus of $\vec{\mathbf c}$. It turns out that $\|
\cdot \|_{\cA^s_{weak}}$ is a quasi--norm and we refer to
\cite{CDD1,DeVore} for further details on the quasi--Banach spaces
$\cA^s_{weak}$. Such spaces can be usually found in literature under the notation $\ell_\tau^w$ (weak--$\ell_\tau$) 
where $\tau =(1/2+s)^{-1} \in (0,2)$, and they are nothing but particular instances of Lorentz sequence spaces. 
Let us only mention that
\begin{equation}
\label{ell2eq}
\| \vec{\mathbf{c}} \|_{\cA^s_{weak}} \sim \sup_{N \in \mathbb{N}} N^s \| \vec{\mathbf{c}} - \vec{\mathbf{c}}_N\|_{\ell_2(\cN)},
\end{equation}
where $\mathbf{c}_N$ is the best $N-$term approximation of
$\vec{\mathbf{c}}$, i.e., the subsequence of $\vec{\mathbf{c}}$
consisting of the $N$ largest coefficients in modulus of
$\vec{\mathbf{c}}$. In particular, \eqref{ell2eq} implies that for
all $\varepsilon>0$ there exists $N_\varepsilon>0$ large enough
such that for all $\vec{\mathbf{c}}$ the best $N_\varepsilon-$term
approximation $\vec{\mathbf{c}}_{N_\varepsilon}$ has the following
properties
\begin{itemize}
\item[(i)] $\|\vec{\mathbf{c}} - \vec{\mathbf{c}}_{N_\varepsilon}\|_{\ell_2{(\cN)}} \leq \varepsilon$;
\item[(ii)] $\# \supp(\vec{\mathbf{c}}_{N_\varepsilon}) \lesssim \varepsilon^{-1/s} \| \vec{\mathbf{c}} \|_{\cA^s_{weak}}^{1/s}$;
\item[(iii)] $\|\vec{\mathbf{c}}_{N_\varepsilon}\|_{\cA^s_{weak}} \leq \| \vec{\mathbf{c}}\|_{\cA^s_{weak}}$.
\end{itemize}
Furthermore, from \eqref{ell2eq} one gets  the following
useful technical estimate
\begin{equation}
\label{finitsupp}
\| \vec{\mathbf{c}} \|_{\cA^{\tilde s}_{weak}} \lesssim (\# \supp( \vec{\mathbf{c}}))^{\tilde s - s} \| \vec{\mathbf{c}} \|_{\cA^{s}_{weak}}
\end{equation}
whenever $0<s< \tilde s$ and $\vec{\mathbf{c}}$ has finite
support. This optimal class of vectors perfectly fits with
complexity estimates for adaptive schemes for elliptic linear
equations \cite{CDD1,CDD2,DFR,S}. For more general nonlinear
problems  a ``weaker'' version of the sparseness
class has been introduced and denoted by $\cA^s_{tree}$ in \cite{CDD3}.

It is not known whether there exist frames for which the solutions
of generic nonlinear equations, particularly for MHD equations,
can have frame coefficients belonging to the sparseness classes
described above. Therefore, we discuss here the requirements,
fulfilled by $\cA_{weak}^s$ and $\cA_{tree}^s$, that a generic
sparseness class should have to ensure quasi--optimality of our
scheme. And, in case the solution belongs to any sparseness class
with such properties, the algorithm will behave as ensured
theoretically. The numerical tests in \cite{U3,U} for turbulent flows
motivate our assumption that the (wavelet) frame coefficients of
the corresponding solutions do belong to some of these generic sparseness
classes (i.e., only few significant wavelet coefficients can be
expected to be relevant in the representation of the solution). Our conceptual approach is also
motivated by the need to simplify the presentation of our
complexity result without going into rather technical details
of the properties of the particular instances $\cA_{weak}^s$ and $\cA_{tree}^s$ of $\cA^s$.
We refer the reader to \cite{CDD3} for more specific details.

For $s>0$ and for a nondecreasing function $T:\mathbb{N} \rightarrow \mathbb{N}$ such that $N \lesssim T(N)$
 we call any space $\cA^s$ a {\it $T$--sparseness class} if $\cA^{\tilde s}_{weak} \subset \cA^s \subseteq \cA^s_{weak}$
 and $\cA^{\tilde s}\subset \cA^s$ for all $\tilde s> s$ and if for all $\vec \uu \in \cA^s$ and for all
 $\varepsilon>0$ there exists a finite vector $\vec \uu_\varepsilon$ with the properties
\begin{itemize}
\item[a)] $\| \vec \uu - \vec \uu_\varepsilon\| \leq \varepsilon$;
\item[b)] $T(\# \supp(\vec \uu_\varepsilon)) \lesssim \varepsilon^{-1/s} \|\vec \uu \|_{\cA^s}^{1/s}$;
\item[c)] $\| \vec \uu_\varepsilon\|_{\cA^s} \lesssim \| \vec \uu \|_{\cA^s}$.
\end{itemize}
In particular we assume that there exists a constant $C_1(s)$ such that
$\|\vec \uu +\vec \vv\|_{\cA^s} \leq C_1(s) (\|\vec \uu\|_{\cA^s} +\|\vec \vv\|_{\cA^s})$.
Of course, $\cA^s_{weak}$ itself is a sparseness class with $T=I$.
Moreover, there exist other T--sparseness classes different from
$\cA^s_{weak}$, for example, the class $\cA^s_{tree}$ defined in
\cite[Formula (6.7)]{CDD3}, that also turns out to be relevant in
our context.

Note that, for $\tilde{\tilde s} > \tilde s > s> 0$, by the
inclusions $\cA^{\tilde{\tilde{s}}}_{weak} \subset \cA^{\tilde s}
\subset \cA^s \subseteq \cA^s_{weak}$ and \eqref{finitsupp} we
have
\begin{equation}
\label{finitsupp2}
\| \vec{\mathbf{c}} \|_{\cA^{\tilde{s}}} \lesssim \| \vec{\mathbf{c}} \|_{\cA^{\tilde{\tilde s}}_{weak}} \lesssim (\#
\supp( \vec{\mathbf{c}}))^{\tilde{ \tilde  s} - s} \| \vec{\mathbf{c}} \|_{\cA^{s}_{weak}} \lesssim  (T(\# \supp(
\vec{\mathbf{c}})))^{\tilde{ \tilde s} - s} \| \vec{\mathbf{c}} \|_{\cA^{s}},
\end{equation}
for all finite vectors $\vec{\mathbf{c}}$.

Now we are ready to formulate our main conceptual requirements. For a fixed $s>0$
\begin{itemize}
\item[(A1)] Let $\theta < 1/3$. We assume that for any $\varepsilon>0$, $\vec{\vv} \in \cA^s$ and
any finitely supported $\vec{\mathbf{w}}$ such that
$$
\| \vec \vv - \vec{\mathbf{w}}\| \leq \theta \varepsilon,
$$
for $\vec{\mathbf{w}}^* = \mathbf{COARSE}[(1-\theta)\varepsilon, \vec{\mathbf{w}}]$ it holds that
$$
    T(\# \supp(\vec{\mathbf{w}}^*)) \lesssim \varepsilon^{-1/s} \| \vec \vv \|_{\cA^s}^{1/s},
$$
and
$$
    \| \vec{\mathbf{w}}^* \|_{\cA^s} \leq C_2(s) \| \vec \vv \|_{\cA^s},
$$
for some costant $C_s(s)>0$.
Moreover, we assume that the number of algebraic operation needed
to compute
$\vec{\mathbf{w}}^*:=\mathbf{COARSE}[\varepsilon,
\vec{\mathbf{w}}]$ for any finite vector $\vec{\mathbf{w}}$ can be
estimated by $\lesssim
T(\#\supp(\vec{\mathbf{w}}))+o(\varepsilon^{-1/s} \|
\vec{\mathbf{w}}\|^{1/s}_{\cA^s})$. See \cite[Proposition 3.2]{S}
and \cite[Proposition 6.3]{CDD3} for examples of procedures
$\mathbf{COARSE}$ with such properties. As it will be clear in the
proof of Theorem \ref{cmpxtm} the use of the procedure
$\mathbf{COARSE}$ is fundamental to ensure that the supports of
the iterates generated by the algorithm can be controlled.
Due to the redundancy of the system used for the discretization, one can expect that the introduction of greedy 
algorithms of matching pursuit type \cite{DMA,GV,MZ,Tem,Tr} can potentially allow even higher rate of compression 
of the active coefficients with respect to the current implementations of $\mathbf{COARSE}$. We postpone this 
investigation to forthcoming work.

\item[(A2)]  The vector $\vec{ \mathbf{w}}_\varepsilon:=\mathbf{APPLY}[\varepsilon,\mathbf{A},\vec{ \mathbf{v}}]$ is such that
$\|\vec{ \mathbf{w}}_\varepsilon\|_{\cA^s} \lesssim \|\vec{ \mathbf{v}}\|_{\cA^s}$, $T(\# \supp(\vec{\mathbf{w}}_\varepsilon)) \lesssim \varepsilon^{-1/s} \| \vec{\mathbf v}\|_{\cA^s}^{1/s}$, and it is computed with a number of
algebraic operations estimable by $\lesssim \varepsilon^{-1/s} \|  \vec{ \mathbf{v}} \|_{\cA^s}^{1/s} + T(\#\supp(\vec{\mathbf{v}}))$.
See \cite[Proposition 3.8]{S} and \cite[Corollary 7.5]{CDD3} for examples of procedures $\mathbf{APPLY}$ with such properties.

\item[(A3)] One of the most crucial procedures of the iterative approximate fixed point scheme is
the realization of $\vec{\mathbf{g}}^{(j)}_i :=
\mathbf{RHS}[\frac{\theta \varepsilon_j}{12 \alpha
K},\vec{\mathbf{\ell}} - \mathbf{A}_1(\vec \vv_i)]$ for each step
$i$ and $j$ of the outer and inner loops, respectively. In
particular, it requires computing efficiently a finite
 approximation to $\mathbf{A}_1(\vec \vv_i)$, where $\ba_1$ is some nonlinear operator and $\vec \vv_i$ is a
 given finite vector. In their
pioneering work \cite{CDD3,CDD4}, Cohen, Dahmen, and DeVore have
found an effective way for solving this problem  for
multiscale and wavelet expansions. 
Note that the efficient evaluation of nonlinear functionals 
(as the ones we consider for the MHD problem) on coefficients as in \cite{CDD3,CDD4} 
is valid as soon as the reference bases are multiscale, i.e., the size of the supports
 of the basis functions decays dyadically for increasing levels, and their corresponding 
 coefficients can be structured in suitable trees. 
Moreover, one needs that such bases characterize Besov spaces by certain norm equivalences.
In Section 7, we present a construction of suitable wavelet frames that enjoy 
the above mentioned properties and, therefore, allow for a straightforward application of the results 
in \cite{CDD3,CDD4}.

We only mention Cohen, Dahmen, and DeVore
results here and refer to \cite{CDD3,CDD4} for technical details (in particular see \cite[Theorem 7.3 and Theorem 7.4]{CDD3}):
there exists a procedure $\mathbf{RHS}$ such that
$\|\vec{\mathbf{g}}^{(j)}_i -(\vec \bl - \mathbf{A}_1(\vec \vv_i))
\|_{\ell_2} \leq \frac{\theta \epsilon_j}{12 \alpha K}$, $T(\#
\supp(\vec{\mathbf{g}}^{(j)}_i)) \lesssim \epsilon_j^{-1/s}
\|\vec{\mathbf{v}}_i\|_{\cA^s}^{1/s}$, $\|\vec{\mathbf{g}}^{(j)}_i\|_{\cA^s} \lesssim 1+ \|\vec{\mathbf{v}}_i\|_{\cA^s}$, and the number of algebraic
operations needed to compute $\vec{\mathbf{g}}^{(j)}_i$ is bounded
by $\lesssim \epsilon_j^{-1/s}
\|\vec{\mathbf{v}}_i\|_{\cA^s}^{1/s} + T(\# \supp(\vec \vv_i))$.
\end{itemize}


In the following the subscript index $i$ refers to the iterations in the outer loop of the fixed point iteration
and the superscript $j$   refers to the inner loop iterations in $\mathbf{SOLVE}$. 
Moreover, $\varepsilon_i$ and $\epsilon_j$ refer to the outer and inner loop tolerances, 
respectively. 
All  estimations below  hold asymptotically for $\varepsilon \rightarrow 0$ (
$\varepsilon$ as in Algorithm \ref{alg2} ).

\begin{tm}
\label{cmpxtm} For $0<s<\tilde s< \tilde{\tilde{s}}$ let $\cA^{\tilde s}$ be a
$T$-sparseness class and  $\vec \uu \in \cA^{s}$,  the
solution of \eqref{abfor2} as in Corollary \ref{discrsolution}.
Assume that 
\begin{itemize}
\item[(i)]  (A1)-(A3) hold for all $s \in (0,\tilde s]$;
\item[(ii)] $\mathbf{P}$ is bounded on $\cA^t$ for all $t \in (0,\tilde s]$;
\item[(iii)] $K>0$ and $0<\theta<1/3$ in Algorithm \ref{alg1} are chosen   so that 
\begin{equation}
\label{Klarge}
 C_1(s) C_2(s)  \| \id-\mathbf P\| (3 \rho^K/\theta)^{\tilde{\tilde{s}}/s-1}<1.
\end{equation}
Here the norm $\| \id-\mathbf P\|$ is the norm of $\id-\mathbf P$ as an operator on $\cA^s$;
\item[(iv)] the constants $L, \varepsilon_0$ and $r^*$ satisfy
\begin{equation}
\label{techniq}
L < \varepsilon_0 < r^* \quad and \quad  \rho^{-K} \frac{L}{\varepsilon_0
-L}+\delta_0\le \frac{1}{2}, \quad for \ some \ \delta_0>0.
\end{equation}
\end{itemize}
 Then for any  $\varepsilon>0$ and  $\delta>0$ such that
$\tilde{\tilde{s}}/\tilde{s}=1+\delta$, the finite vector $\vec
\uu_\varepsilon:= \mathbf{FIXPT}[\varepsilon, \ba,\bA1,\vec{\bl}]$
satisfies
\begin{itemize}
\item[a)] $\| \vec \uu- \mathbf P \vec \uu_\varepsilon \|_{\ell_2(\cN)} \leq \varepsilon$;
\item[b)] $\#\supp(\vec \uu_\varepsilon) \lesssim \varepsilon^{-(1+\delta)/s} \| \vec \uu\|_{\cA^s}^{(1+\delta)/s}$;
\item[c)] the number of algebraic operations needed to compute $\vec \uu_\varepsilon$ is $\lesssim \varepsilon^{-(1+\delta)/s} \|\vec \uu \|_{\cA^s}^{(1+\delta)/s}$.
\end{itemize}
\end{tm}

\begin{proof}
For the proof of part a) see Theorem \ref{converg}. Next, we show part b). 
Assume that $\{\vec \vv_i\}_{i \in \NN_0}$
is the sequence of vectors generated in Algorithm \ref{alg2}.
We want to show that $T(\# \supp(\vec \vv_i)) \lesssim
\varepsilon_i^{-(1+\delta)/s} \|\vec \uu\|_{\cA^s}^{(1+\delta)/s}$
for  $i$ large enough. Since $\mathbf{P}$ is bounded on $\cA^t$ for
all $t \in (0,\tilde s]$, it is also bounded  on $\cA^s$.
Therefore $\mathbf{P}  \vec \uu \in \cA^s$. Then for
$\epsilon_j>0$ there exists a finite vector $(\mathbf{P}  \vec
\uu)_{\epsilon_j}$ such that $\| \mathbf{P}  \vec \uu -
(\mathbf{P}  \vec \uu)_{\epsilon_j}\|_{\ell_2(\cN)} \leq
\frac{\theta}{6}\epsilon_j$ and $\# \supp((\mathbf{P} \vec
\uu)_{\epsilon_j}) \lesssim T(\# \supp((\mathbf{P}  \vec
\uu)_{\epsilon_j})) \lesssim \epsilon_j^{-1/s} \| \mathbf{P}  \vec
\uu\|_{\cA^s}^{1/s} \lesssim \epsilon_j^{-1/s} \|\vec
\uu\|_{\cA^s}^{1/s}$. Therefore, by \eqref{finitsupp2} we have
\begin{equation}
\label{finitsupp3}
\epsilon_j^{\tilde{\tilde{s}}/s-1}\|(\mathbf{P}  \vec \uu)_{\epsilon_j}\|_{\cA^{\tilde s}}
\lesssim \|\vec \uu\|_{\cA^s}^{\tilde{\tilde{s}}/s-1}\|(\mathbf{P}  \vec \uu)_{\epsilon_j}\|_{\cA^{s}}
\lesssim \|\vec \uu\|_{\cA^s}^{\tilde{\tilde{s}}/s-1}\|\mathbf{P}  \vec \uu\|_{\cA^{s}} \lesssim
\|\vec \uu\|_{\cA^s}^{\tilde{\tilde{s}}/s-1}\|\vec \uu\|_{\cA^{s}} = 
\|\vec \uu\|_{\cA^s}^{\tilde{\tilde{s}}/s}
\end{equation}
for any $\tilde{\tilde{s}} > \tilde{s}>s>0$. From
\eqref{SOLVEhelp} we get that
\begin{equation}
\label{SOLVEhelp2}
\bigl\|\mathbf P \vec{\mathbf v}_i^{ex}-(\id-\mathbf P)\vec{\mathbf{v}}_i^{(j-1)}-
\vec{\mathbf{v}}_i^{(j,K)}\bigr\|_{\ell_2(\cN)}\le\frac{2\theta\epsilon_j}{3}
\end{equation}
with $\vec{\mathbf v}_i^{ex}:= \mathbf{H}
(\vec{\mathbf{v}}_{i-1})$. Due to \eqref{techniq} and
\eqref{iteration1}, we obtain for any $i$ large enough and some $\delta_0 >0$ that
\begin{eqnarray*}
 \| \mathbf P\vec{\mathbf u} -\mathbf P\vec{\mathbf v}_i^{ex}\|_{\ell_2(\cN)}  &\leq &\| \mathbf{H} ( \vec{\mathbf u}) - \mathbf{H}(\vec{\mathbf{v}}_{i-1})\|_{\ell_2(\cN)} \\
&\leq& L \|\mathbf{P} \vec{\mathbf u} - \mathbf{P}  \vec{\mathbf{v}}_{i-1} \|_{\ell_2(\cN)}  \\
&\leq & L \L \frac{\varepsilon_0^{i-1} ( \varepsilon_0 - L (L/\varepsilon_0)^{i-2})}{\varepsilon_0 - L} + L^{i-2} r^* \R\\
&=& \frac{L}{\varepsilon_0}\L \frac{\varepsilon_0^{i} ( \varepsilon_0 - L (L/\varepsilon_0)^{i-2})}{\varepsilon_0 - L} +
 L^{i-2} \varepsilon_0 r^* \R \\
&=& \varepsilon_0^{i}  \frac{L}{\varepsilon_0}\L \frac{( \varepsilon_0 - L (L/\varepsilon_0)^{i-2})}{\varepsilon_0 - L}
+ \L \frac{L}{\varepsilon_0}\R^{i-2} \frac{r^*}{\varepsilon_0} \R \\
&\leq& \varepsilon_0^{i} \L \frac{L}{\varepsilon_0
-L}+\delta_0 \R.
\end{eqnarray*}
Due to the stopping criterion of Algorithm \ref{alg1} we have for
all $i$ that for the last $j-$th iteration $\varepsilon_0^i \leq 
\frac{\theta }{3 \rho^K} \epsilon_j$. Therefore, for all $i$ large enough, due to
\eqref{techniq} we get
$$
  \| \mathbf P\vec{\mathbf u} -\mathbf P\vec{\mathbf v}_i^{ex}\|_{\ell_2(\cN)} \le \frac{\theta}{6} \epsilon_j,
$$
and
\begin{equation}
\label{SOLVEhelp3} \bigl\|\mathbf P\vec{\mathbf u} -(\id-\mathbf
P)\vec{\mathbf{v}}_i^{(j-1)}-\vec{\mathbf{v}}_i^{(j,K)}\bigr\|_{\ell_2(\cN)}
\le \frac{\theta \epsilon_j}{6}+\frac{2 \theta \epsilon_j}{3},
\end{equation}
which implies that
 \begin{equation}
\label{SOLVEhelp5}
\bigl\|(\mathbf P\vec{\mathbf u})_{\epsilon_j} -(\id-\mathbf P)\vec{\mathbf{v}}_i^{(j-1)}-
\vec{\mathbf{v}}_i^{(j,K)}\bigr\|_{\ell_2(\cN)}\le \theta \epsilon_j.
\end{equation}
Due to  (A1), \eqref{SOLVEhelp5}, and   $\|\cdot \|_{\cA^{\tilde s}}$ being a quasi-norm, it follows that
$\vec{\mathbf{v}}_i^{(j)}:=\mathbf{COARSE}[(1-\theta)\epsilon_j,\vec{\mathbf{v}}_i^{(j,K)}]$, for $i$ large enough,
satisfies
\begin{eqnarray*}
\| \vec{\mathbf{v}}_i^{(j)} \|_{\cA^{\tilde s}} &\leq& C_2(s)\|(\mathbf P\vec{\mathbf u})_{\epsilon_j} -(\id-\mathbf P)\vec{\mathbf{v}}_i^{(j-1)}\|_{\cA^{\tilde s}} \\
&\leq& C_1(s) C_2(s)  \L \|(\mathbf P\vec{\mathbf u})_{\epsilon_j}\|_{\cA^{\tilde s}} + \| \id-\mathbf P\| \| \vec{\mathbf{v}}_i^{(j-1)}\|_{\cA^{\tilde s}}\R,
\end{eqnarray*}
so by \eqref{finitsupp3} and $\epsilon_j = 3 \rho^K/\theta \epsilon_{j-1}$ (see Algorithm \ref{alg1}),
\begin{eqnarray*}
\L\epsilon_j^{\tilde{\tilde{s}}/s-1} \| \vec{\mathbf{v}}_i^{(j)}\|_{\cA^{\tilde s}} \R
&\leq&  C' \|\vec{\mathbf u}\|_{\cA^{s}}^{\tilde{\tilde{s}}/s} + C_1(s) C_2(s)  \| \id-\mathbf P\| (3 \rho^K/\theta)^{\tilde{\tilde{s}}/s-1}\L\epsilon_{j-1}^{\tilde{\tilde{s}}/s-1}\| \vec{\mathbf{v}}_i^{(j-1)}\|_{\cA^{\tilde s}} \R.
\end{eqnarray*}
We can conclude that for $K>0$ large enough, and by the assumption \eqref{Klarge}, the solutions of the homogeneous part of this recursion converge to zero, and so
\begin{equation}
\label{boundedseq}
\epsilon_j^{\tilde{\tilde{s}}/s-1} \| \vec{\mathbf{v}}_i^{(j)}\|_{\cA^{\tilde s}} \lesssim \|\vec{\mathbf u}\|_{\cA^{s}}^{\tilde{\tilde{s}}/s},
\end{equation}
uniformly with respect to $j$.  And, due to (A1) we also have
\begin{eqnarray*}
\# \supp( \vec{\mathbf{v}}_i^{(j)}) &\lesssim& T( \# \supp( \vec{\mathbf{v}}_i^{(j)})) \\
&\lesssim &\epsilon_j^{-1/\tilde s} \| (\mathbf P\vec{\mathbf u})_{\epsilon_j} -(\id-\mathbf P)\vec{\mathbf{v}}_i^{(j-1)}\|^{1/{\tilde s}}_{\cA^{\tilde s}}\\
&\lesssim& \epsilon_j^{-\frac{\tilde{\tilde s}/\tilde s}{s}} \L
\epsilon_j^{\frac{\tilde{\tilde s}}{s}-1} \left [ \| (\mathbf
P\vec{\mathbf u})_{\epsilon_j}\|_{\cA^{\tilde s}} + \|
\id-\mathbf P \| \| \vec{\mathbf{v}}_i^{(j-1)}\|_{\cA^{\tilde s}}
\right ]\R^{1/\tilde s} .
\end{eqnarray*}
Therefore, by  \eqref{finitsupp3} and \eqref{boundedseq} we get
with $\tilde{\tilde{s}}/\tilde{s}=1+\delta$ that
\begin{equation}\label{boundedsup}
 \# \supp( \vec{\mathbf{v}}_i^{(j)}) \lesssim T( \# \supp( \vec{\mathbf{v}}_i^{(j)})) \lesssim \epsilon_j^{-\frac{\tilde{\tilde s}/\tilde s}{s}}
 \| \vec{\mathbf u}\|_{\cA^{s}}^{\frac{\tilde{\tilde s}/\tilde s}{s}}
 = \epsilon_j^{-\frac{1+\delta}{s}} \| \vec{\mathbf u}\|_{\cA^{s}}^{\frac{1+\delta}{s}}.
\end{equation}
Next, recall that $0<s<\tilde{\tilde{s}}$ and by Algorithm
\ref{alg1} we have $\varepsilon_i \lesssim \epsilon_j$ for all
$j$. Then \eqref{boundedseq} implies that for all $j$  and $i$ large enough we have $ \varepsilon_i^{\tilde{\tilde{s}}/s-1} \|
\vec{\mathbf{v}}^{(j)}_i\|_{\cA^{\tilde s}} \lesssim
\|\vec{\mathbf u}\|_{\cA^{s}}^{\tilde{\tilde{s}}/s}$. In
particular, for $i$ large enough,
\begin{equation}
\label{boundedseq2} \varepsilon_i^{\tilde{\tilde{s}}/s-1} \|
\vec{\mathbf{v}}_i\|_{\cA^{\tilde s}} \lesssim \|\vec{\mathbf
u}\|_{\cA^{s}}^{\tilde{\tilde{s}}/s}.
\end{equation}
By the same argument, \eqref{boundedsup} yields for $i$ large enough
\begin{equation}\label{boundedsup2}
 \# \supp( \vec{\mathbf{v}}_i) \lesssim T( \# \supp( \vec{\mathbf{v}}_i)) \lesssim
 \varepsilon_i^{-\frac{1+\delta}{s}} \| \vec{\mathbf u}\|_{\cA^{s}}^{\frac{1+\delta}{s}}.
\end{equation}
Note that the stopping criterion in Algorithm \ref{alg2} implies
that for large enough $i$ we have $\varepsilon \lesssim
\varepsilon_i$. Therefore, from  \eqref{boundedsup2} we also get
that
\begin{equation}\label{boundedsup3}
 \# \supp( \vec{\mathbf{u}}_\varepsilon) \lesssim T( \# \supp( \vec{\mathbf{u}}_\varepsilon)) \lesssim
 \varepsilon^{-\frac{1+\delta}{s}} \| \vec{\mathbf u}\|_{\cA^{s}}^{\frac{1+\delta}{s}}.
\end{equation}

To prove part c) it is sufficient to show that the number of
algebraic operations needed for each iterations, i.e., for the
computation of $\vec{\mathbf{v}}_i:=
\mathbf{SOLVE}[\varepsilon_{i},\ba,\vec{\bl}-\bA1(\vec
\vv_{i-1})]$, is  $\lesssim \varepsilon_i^{-(1+\delta)/s} \|\vec
\uu\|_{\cA^s}^{(1+\delta)/s}$ for all $i$ large enough. 
Note that for small  $i$  the number of
algebraic operations needed  is bounded by $C \varepsilon^{-(1+\delta)/s} \|\uu\|_{\cA^s}^{(1+\delta)/s}$ 
for $\varepsilon \rightarrow 0$.

Due to $L<1$ and by
\eqref{estimation2}, there exists $N \in \NN$ such that $L^N r^*
\leq \varepsilon$. Therefore, the $N-$th, exit iteration,
satisfies $N:=\left [\frac{\log_{\varepsilon_0^{-1}}(\varepsilon/r^*)}{\log_{\varepsilon_0^{-1}}(L)} \right]
<\left [-\log_{\varepsilon_0^{-1}}(\varepsilon)-\frac{\log_{\varepsilon_0^{-1}}(r^*)}{\log_{\varepsilon_0^{-1}}(L)} \right]$.
This implies that we can estimate the total number of
operations by
\begin{eqnarray*}
\lesssim \sum_{i=0}^{\left [-\log_{\varepsilon_0^{-1}}(\varepsilon)-\frac{\log_{\varepsilon_0^{-1}}(r^*)}{\log_{\varepsilon_0^{-1}}(L)} \right]} \varepsilon_i^{-(1+\delta)/s} \|\vec \uu\|_{\cA^s}^{(1+\delta)/s}
&=& \|\uu\|_{\cA^s}^{(1+\delta)/s} \sum_{i=0}^{\left [-\log_{\varepsilon_0^{-1}}(\varepsilon)-\frac{\log_{\varepsilon_0^{-1}}(r^*)}{\log_{\varepsilon_0^{-1}}(L)}\right]} \L\varepsilon_0^{-(1+\delta)/s}\R^i \\
&=& \|\uu\|_{\cA^s}^{(1+\delta)/s} \frac{1- \L\varepsilon_0^{-(1+\delta)/s}\R^{[-\log_{\varepsilon_0^{-1}}(\varepsilon)-\frac{\log_{\varepsilon_0^{-1}}(r^*)}{\log_{\varepsilon_0^{-1}}(L)}]+1}}{1-\varepsilon_0^{-(1+\delta)/s}}\\
&\lesssim & \varepsilon^{-(1+\delta)/s} \|\uu\|_{\cA^s}^{(1+\delta)/s}.
\end{eqnarray*}
The last inequality is due to the fact that here we consider the asymptotic behavior for 
$\varepsilon \rightarrow 0$. 

Therefore to conclude the
proof, it is sufficient to estimate the complexity of
$\mathbf{SOLVE}$. To do so one can follow the argument used in the
proof of \cite[Theorem 3.12 (II)]{S} replacing \cite[(3.19)]{S} by
\eqref{boundedseq} and using the assumptions (A1)-(A3) when
relevant.  
Due to assumption (A3) one has $T(\#\supp(\vec{\mathbf g}^{(j)}_{i-1})) \lesssim \epsilon_j^{-1/\tilde s} \|\vec{\mathbf{v}}_{i-1}\|_{\cA^{\tilde s}}^{1/\tilde s}$ and $\| \vec{\mathbf g}^{(j)}_{i-1} \|_{\cA^{\tilde s}} \lesssim (1 + \| \vec{\mathbf{v}}_{i-1}\|_{\cA^{\tilde s}})$.
Therefore by assumption (A2) one has $T(\#\supp(\vec{\mathbf f}^{(j)}_{i-1})) \lesssim \epsilon_j^{-1/\tilde s} \|\vec{\mathbf{v}}_{i-1}\|_{\cA^{\tilde s}}^{1/\tilde s}$ and $\| \vec{\mathbf f}^{(j)}_{i-1} \|_{\cA^{\tilde s}} \lesssim (1 + \| \vec{\mathbf{v}}_{i-1}\|_{\cA^{\tilde s}})$.
By induction assumption we have that 
$\| \vec{\mathbf v}^{(j-1)}_i \|_{\cA^{\tilde s}} \lesssim (1 + \| \vec{\mathbf{v}}_{i-1}\|_{\cA^{\tilde s}})$ 
and $T(\# \supp(\vec{\mathbf v}^{(j-1)}_i)) \lesssim \epsilon_j^{-1/\tilde s} \| \vec{\mathbf{v}}_{i-1}\|_{\cA^{\tilde s}}^{1/\tilde s}$, that are trivially valid for $j=1$. Therefore, again by (A2) one has
\begin{equation} \label{aux1}
\| \vec{\mathbf v}^{(j,k)}_i \|_{\cA^{\tilde s}} \lesssim (1 + \| \vec{\mathbf{v}}_{i-1}\|_{\cA^{\tilde s}}),
\end{equation}
 and
\begin{equation} \label{aux2}
T(\# \supp(\vec{\mathbf v}^{(j,k)}_i)) \lesssim \epsilon_j^{-1/\tilde s} \| \vec{\mathbf{v}}_{i-1}\|_{\cA^{\tilde s}}^{1/\tilde s},
\end{equation}
for all $0 \leq k \leq K$. Therefore, by (A2) and (A3), the number of algebraic operations required to 
compute $\vec{\mathbf v}^{(j,K)}_i$ is $\lesssim \epsilon_j^{-1/\tilde s} 
(1 + \| \vec{\mathbf{v}}_{i-1}\|_{\cA^{\tilde s}})^{1/\tilde s} + T(\# \supp(\vec{\mathbf v}_{i-1}))$. 
Finally by assumption (A1) the application of $\vec{\mathbf v}_i^{(j)}:=
\mathbf{COARSE}[(1-\theta)\epsilon_j,\vec{\mathbf v}_i^{(j,K)}]$ costs $\lesssim 
T(\# \supp(\vec{\mathbf v}_i^{(j,K)})) + o(\epsilon_j^{-1/\tilde s} \| \vec{\mathbf v}_i^{(j,K)}
\|_{\cA^{\tilde s}}^{1/\tilde s}) \lesssim  \epsilon_j^{-1/\tilde s} (1+\| \vec{\mathbf{v}}_{i-1}
\|_{\cA^{\tilde s}})^{1/\tilde s}$. Thus, by (A1), \eqref{aux1} and \eqref{aux2} the induction
hypothesis holds, i.e. $\| \vec{\mathbf v}^{(j)}_i \|_{\cA^{\tilde s}} 
\lesssim (1 + \| \vec{\mathbf{v}}_{i-1}\|_{\cA^{\tilde s}})$ and $T(\# \supp(\vec{\mathbf v}^{(j)}_i)) 
\lesssim \epsilon_j^{-1/\tilde s} \| \vec{\mathbf{v}}_{i-1}\|_{\cA^{\tilde s}}^{1/\tilde s}$ for any $j \in \NN_0$.
This implies by induction that computing 
$\vec{\mathbf v}^{(j)}_i$ from $\vec{\mathbf v}^{(j-1)}_i$ takes a number of operations 
$\lesssim  \epsilon_j^{-1/\tilde s} (1 + \| \vec{\mathbf{v}}_{i-1}\|_{\cA^{\tilde s}})^{1/\tilde s} + 
T(\# \supp(\vec{\mathbf v}_{i-1}))$.
Since $\epsilon_j$ decreases geometrically and by formulas \eqref{boundedseq2} and \eqref{boundedsup2}, 
one can estimate, as done before, the cost of the computation of $\vec{\mathbf{v}}_i:=
\mathbf{SOLVE}[\varepsilon_{i},\ba,\vec{\bl}-\bA1(\vec \vv_{i-1})]$ as a multiple of 
\begin{eqnarray*}
&& \varepsilon_i^{-1/\tilde s}   (1 + \| \vec{\mathbf{v}}_{i-1}\|_{\cA^{\tilde s}})^{1/\tilde s} + \log(\varepsilon_i) T(\# \supp(\vec{\mathbf v}_{i-1})) \\
&\lesssim& \varepsilon_i^{-1/\tilde s}   (1 + \| \vec{\mathbf{v}}_{i-1}\|_{\cA^{\tilde s}})^{1/\tilde s} + \log(\varepsilon_i)  \varepsilon_{i-1}^{-\frac{1+\delta}{s}} \| \vec{\mathbf u}\|_{\cA^{s}}^{\frac{1+\delta}{s}}\\
&\lesssim& \varepsilon_i^{-\frac{\tilde{\tilde s}/\tilde s}{s}}   \left (\varepsilon_i^{\frac{\tilde{\tilde s}}{s}-1} (1 + \| \vec{\mathbf{v}}_{i-1}\|_{\cA^{\tilde s}}) \right)^{1/\tilde s} + \log(\varepsilon_i)  \varepsilon_{i-1}^{-\frac{1+\delta}{s}} \| \vec{\mathbf u}\|_{\cA^{s}}^{\frac{1+\delta}{s}}\\
&\lesssim& \varepsilon_i^{-\frac{\tilde{\tilde s}/\tilde s}{s}}   \| \vec{\mathbf u}\|_{\cA^{s}}^{\frac{\tilde{\tilde s}/\tilde s}{s}} + \log(\varepsilon_i)  \varepsilon_{i-1}^{-\frac{1+\delta}{s}} \| \vec{\mathbf u}\|_{\cA^{s}}^{\frac{1+\delta}{s}} \\
&\lesssim &  \varepsilon_{i}^{-\frac{1+\delta}{s}} \| \vec{\mathbf u}\|_{\cA^{s}}^{\frac{1+\delta}{s}}.
\end{eqnarray*}
In the last inequality we could ignore the $\log$ factor because of the arbitrary choice of $\delta>0$, small as one wants. This concludes the proof.
\end{proof}

\begin{rems} 1. Theorem \ref{cmpxtm} ensures that our scheme is arbitrarily close 
to   optimality, i.e., it is quasi--optimal;\\
2.
On the one hand, it is a very difficult theoretical problem to prove that $\mathbf{P}$ is
bounded on $\cA^t$ for all $t \in (0,\tilde{\tilde s})$ for wavelet frames.
On the other hand, this condition has been verified numerically in \cite{W,DFRSW} in case of wavelet frame discretization, 
by observing the optimal convergence of ${\bf SOLVE}$.  According to \cite[Remark 3.13]{S}, the  boundedness of
 $\mathbf{P}$ on $\cA^t$ for all $t \in (0,\tilde{\tilde s})$ is (almost) a necessary requirement for 
the scheme to behave optimally.

There exists frames, for example time--frequency 
localized Gabor frames (and more generally all intrisically polynomially localized frames 
\cite{FG,DFR}), for which the boundedness of the corresponding $\mathbf{P}$ has been proven rigorously, see 
\cite[Theorem 7.1 in Section 7]{DFR}. Therefore, for certain operator equations 
(for example, certain pseudodifferential equations appearing, e.g., in wireless 
communication modeling \cite{G}) the optimal application of $\mathbf{SOLVE}$ based on Gabor frame 
discretizations is justified theoretically.

One can avoid requiring the boundedness of $\PP$ and obtain a fully
optimal scheme replacing $\mathbf{SOLVE}$ by its
modified version $\mathbf{modSOLVE}$ (see \cite{S}). One assumes then, without loss of generality, 
that the number of algebraic
operations needed to compute $\vec{\mathbf{g}}^{(j)}_i$ in   (A3) is bounded  
by $\lesssim \epsilon_j^{-1/s} \|\vec{\mathbf{v}}_i\|_{\cA^s}^{1/s}$. 
The procedure $\mathbf{modSOLVE}$ requires the definition of an implementable alternative projector 
$\tilde \PP$ onto $\Ran(\mathbf{A})$ and it  is possible if a suitable wavelet frame 
expansion is constructed. In Section 7, we present such a construction  that 
allows for a straightworward definition of the admissible $\tilde \PP$ following the recipe in \cite{S}.
\\
If $\mathcal{F}$ is a Riesz basis, then the scheme is certainly optimal and our 
result confirms the one in \cite[Theorem 7.5]{CDD3}. 

3. Due to Theorem \ref{discrprob}, \eqref{techniq}  holds
 if the physical parameters of the MHD problem allow for choosing $L$,
$\varepsilon_0$ and $r^*$ such that
$$
 0 < L=: r^* \gamma < \varepsilon_0 < r^* < \gamma^{-1} 
$$
for $\gamma=\gamma(\mathcal{F},A,A_1):=  2 \|\ba^{-1}|_{\Ran(\ba)}\| \|A_1\| \|F\|^3$. 
In particular, if $0<\gamma<1$ is small enough then  \eqref{techniq} is satisfied.  We   show next that, for 
the particular MHD equations and in case $\mathcal{F}$ is a Riesz basis, the 
  dependence of $\gamma$ on the viscosity $\eta$ and the electric 
  resistivity $\sigma^{-1}$ can be expressed explicitly. Note that
\begin{eqnarray*}
\langle \mathbf{A} \vec{\mathbf u},\vec{\mathbf u}\rangle &=& \langle ( D_{\bV}^*)^{-1} FA F^* D_{\bV}^{-1} \vec{\mathbf u},\vec{\mathbf u}\rangle\\
&=& \langle A F^* D_{\bV}^{-1} \vec{\mathbf u},F^* D_{\bV}^{-1} \vec{\mathbf u}\rangle\\
&\geq & \alpha \|  F^* D_{\bV}^{-1} \vec{\mathbf u}\|_{\bV}\\
&=& \alpha \| \sum_n (  D_{\bV}^{-1} \vec{\mathbf u})_n f_n \|_{\bV} \sim  \alpha \|  D_{\bV}^{-1} \vec{\mathbf u} \|_{\bV_d} = \alpha \| \vec{\mathbf u} \|_{\ell_2(\cN)}.
\end{eqnarray*}
 Therefore, if $\mathcal{F}$ is a Riesz basis, then $\Ran(\mathbf{A}) = \ell_2(\mathcal{N})$ and 
\begin{equation}
\label{alphaestim}
\|\ba|_{\Ran(\ba)}^{-1}\| = \|\ba^{-1}\| \lesssim \alpha^{-1},
\end{equation}
where from \eqref{coercive} we have $\alpha := \Big(\alpha_0-2\|a_1\| \ 
\|(\uu_0, \JJ_0)\|_{\XX_{(\uu,\JJ)}}\Big)$ and $\alpha_0=c(\Omega)\min\{\eta, 
\sigma^{-1}\}$. This implies that if $\eta, \sigma^{-1}$ are large enough then $\gamma<1$
 can be made sufficiently small. The norm equivalence $\| \sum_{n \in \cN} (  D_{\bV}^{-1} \vec{\mathbf u})_n f_n 
\|_{\bV} \sim \|  D_{\bV}^{-1} \vec{\mathbf u} \|_{\bV_d}$ is valid only 
if $\cF$ is a Riesz basis and it does not hold for pure frames. It is 
 possible to get an  estimate  similar to \eqref{alphaestim} under the additional 
assumption on the frame $\mathcal{F}$ that $D_{\bV}^{-1} \Ran(\mathbf{A}) \subseteq 
\Ran(\mathbf{F})$, which is to be verified for any particular frame under 
consideration.
\end{rems}

\section{Construction of aggregated wavelet frames }
In this section we present the construction of suitable multiscale frames on bounded domains for the MHD problem
described in Sections 2--3. They are in fact Gelfand frames for
$(\bV,\cH,\bV')$, which ensures the correct application of
Algorithms \ref{alg1} and \ref{alg2}. In particular, we show that
for the special case  $\Gamma=\Gamma_1=\Sigma_1$ 
we can use the wavelet bases constructed using \cite[Chapter 2, Definition
5]{U} from a biorthogonal system on $L_2(\Omega)$ satisfying
\cite[Chapter 2, Assumption 3]{U}, where
$\Omega$ is some open and bounded domain, see for example \cite{DKU,DS1,DS2}.

We start by studying the properties of the functions in $\bV$. Using  
$$
 \bV:=\{(\vv,\KK) \in \XX_{(\vv,\KK)}: \ b( (\vv,\KK),(q,\psi))=0 \ \text{for all} \ (q,\psi) \in
 M_{(q,\psi)}\},
$$
and setting $\psi=0$ in the definition of the form $b$ we
obtain
$$
 \int_\Omega (\nabla \cdot \vv) q=0 \ \ \text{for all} \ \ q \in \dot{L}_2(\Omega).
$$
By \cite[Lemma 1.a)]{CMS}  we get that
$\nabla \cdot \vv=0$. Thus, $\vv \in \bV(\text{div}; \Omega)$
satisfying $\vv|_\Gamma=0$, where
\begin{eqnarray*}
 \bV(\text{div}; \Omega):=\{\vv \in \LL_2(\Omega): \ \nabla \cdot \vv
 =0\}.
\end{eqnarray*}
Therefore, $\vv$ is in $\bV(\text{div};\Omega) \cap
\bH^1_0(\Omega)$. A detailed discussion of construction of a
wavelet basis for such a space is given, for example, in
\cite{U2,U1,U}.

Note that the application of $\nabla$ to $\psi$ yields the same
results regardless if $\psi \in \dot{H}^1(\Omega)$ or $\psi \in
H^1(\Omega) \setminus \RR$. Thus, substituting $q=0$ into the
definition of the form $b$ we obtain
\begin{equation} \label{sec10:aux1}
 \int_\Omega \KK (\nabla \psi) =0 \ \ \text{for all} \ \ \psi \in H^1(\Omega)\setminus \RR.
\end{equation}
Corollary 3.4. in \cite{GR}  yields, for simply--connected
$\Omega$, that $\KK= \nabla \times \xi$ with $\xi \in
\bH(\text{curl}; \Omega)$, where $
 \bH(\text{curl}; \Omega):=\{\xi \in \LL_2, \ \nabla \times \xi \in \LL_2\}$.
There is also another way to describe the space, to which 
$\KK$ belongs.

\begin{prop} \label{prop:KK} Let $\KK $ and $\psi$ be in $\LL_2(\Omega)$ and
$\dot{H}^1(\Omega)$, respectively. Then
$$
 \int_\Omega \KK (\nabla \psi) =0 \ \ \text{for all} \ \ \psi \in \dot{H}^1(\Omega)
$$
if and only if $\nabla \cdot \KK=0$ and $\KK \cdot \nn|_\Gamma=0$.
\end{prop}
\begin{proof}
''$\Longrightarrow$`` Integration by parts yields
\begin{equation} \label{aux}
 \int_\Omega (\nabla \cdot \KK) \psi=-\int_\Omega \KK (\nabla
 \psi)+\int_\Gamma
 \KK \cdot \nn \psi \ \ \text{for all} \ \psi \in H^1(\Omega) \setminus \RR.
\end{equation} And, by the assumption together with (\ref{sec10:aux1})  we
get that \eqref{aux} reduces to
$$
 \int_\Omega (\nabla \cdot \KK) \psi=0 \ \ \text{for all} \ \psi \in H_0^1(\Omega)
  \subset H^1(\Omega)\setminus \RR .
$$
This implies that $\nabla \cdot \KK=0$ due to $H^{-1}(\Omega)$
being a dual for $H_0^1(\Omega)$. Thus, \eqref{aux} becomes
$$
 \int_\Gamma \KK \cdot \nn \psi=0 \ \ \text{for all} \ \psi \in H^1(\Omega) \setminus \RR,
$$
which implies $\KK \cdot \nn|_\Gamma=0$ as $\KK \cdot \nn|_\Gamma$
is in $H^{1/2}(\Gamma)'$, the dual of $H^{1/2}(\Gamma)$, and $\psi
\in H^{1/2}(\Gamma)$. The claim follows. ''$\Longleftarrow$``
follows directly from \eqref{aux}.
\end{proof}

In the notation of \cite[ (2.26)]{U}, by Proposition
\ref{prop:KK} $\KK \in \bV_0(\text{div};
\Omega):=\bH_0(\text{div}; \Omega) \cap \bV(\text{div}; \Omega)$, where
by \cite[Theorem 2.6]{GR}
$$ \bH_0(\text{div}; \Omega)=\{\KK \in \LL_2(\Omega): \ \nabla
\cdot \KK \in \LL_2(\Omega), \ \KK  \cdot \nn|_\Gamma=0\}.
$$ 
Thus, for example, in case $\Omega=(0,1)^n$, the existence of a Riesz basis for
$\bV_0(\text{div};\Omega)$ follows from \cite[Theorem 10 and Proposition
5, p. 100]{U}. Therefore, we are guaranteed that a wavelet frame
needed for discretizing {\bf Problem 3} exists, at least, for any domain which is an 
affine image of $(0,1)^n$ \cite[p. 103]{U}. Restricting our problem to 2D, we could even work with 
divergence--free wavelets on domains $\Omega \subset \mathbb{R}^2$ which are conformal 
images of $(0,1)^2$ \cite[p. 104]{U}. 
\\

We discuss next the methods for
obtaining   pure frames for $\bV$ on more general domains. The following two approaches are
of interest: 
\begin{itemize}
\item Modify the construction of pure
frames for $\bH^s(\Omega)$ given in \cite{DFR} and done using the
ODD (Overlapping Domain Decomposition) technique. Using ODD assume 
that $\{\Omega_i\}_{i=1}^M$, $M \in \NN$, are overlapping subdomains 
such that $ \Omega=\cup_{i=1}^M\Omega_i$. Such domains 
could be assumed to be affine images of the reference domain $\Box:=(0,1)^n$ (see Figure 4).
\begin{figure}[ht]
\hbox to \hsize {\hfill \epsfig{file=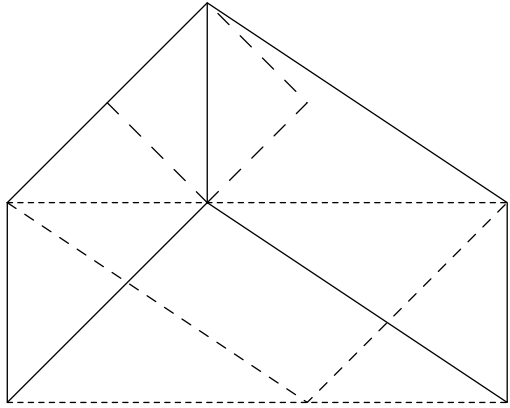,width=7cm} \hfill}
\caption{Example of an Overlapping Domain Decomposition of a polygonal domain in 2D by means of 
patches which are affine images of $(0,1)^n$.}
\end{figure}
Moreover,   assume that for each $1 \le i \le M$ a divergence--free wavelet basis 
${\bf \Psi}_i:=\{{\psi}_{j,k}^i\}_{j \geq -1, k \in \mathcal{J}_j^i}$ is given for 
$\mathbf{V}(\text{div};\Omega_i)\cap \mathbf{H}^1_0(\Omega_i)$. Show that 
${\bf \Psi}:=\cup_{i=1}^N {\bf \Psi}_i = \{ {\psi}_{j,k}^i\}_{j \geq -1, k \in 
\mathcal{J}_j^i, i=1,...,M}$ is   a Gelfand frame for $\mathbf{V}(\text{div};
\Omega)\cap \mathbf{H}^1_0(\Omega)$.
\item Given a wavelet Gelfand frame for $H^1(\Omega)$, possibly constructed 
by ODD, follow a similar strategy as in \cite{U2,U1,U} and construct a divergence-free 
Gelfand frame for $\mathbf{V}(\text{div};\Omega)\cap\bH^1_0(\Omega)$.
\end{itemize}

These systems, called {\it aggregated divergence--free wavelet frames}, 
allow us to avoid dealing with interfacing patches used in disjoint
 domain decompositions (DDD) (see \cite{DS2,U}). DDD are usually rather complicated to 
 implement and can yield ill--conditioned   systems (see \cite[p. 104, sec. More General Domains]{U}).\\

We continue by presenting some preliminary  results  on  constructions of the desired pure wavelet frames
described in the first approach above. Assume  that
$\cC:=\{\Omega_i\}_{i=1}^M$ is an overlapping, relatively compact covering of $\Omega$ such that
\begin{itemize}
\item [(C1)] there exist affine (or even conformal if $n = 2$) maps
  $ \kappa_i:\Box\rightarrow\Omega_i$, for all $i=1,\ldots,M$.
\end{itemize}

The set of admissible domains $\Omega$ is restricted by raising  condition (C1), e.g.,
the boundary of $\Omega$ has to be piecewise smooth.
Nevertheless, the particularly attractive case of polyhedral domains is still covered. Thus,
the above assumptions on the parametrizations $\kappa_i$ are by no means  principal restrictions on the
shape of $\Omega$.

We assume that the wavelets we use below have sufficient regularity 
and number of vanishing moments. We consider a template vector--valued wavelet divergence--free basis
${\bf \Psi}^\Box=\{\psi^\Box_{j,k}\}_{j \geq j_0, k \in \mathcal{J}_j^\Box}$ for $\mathbf{V}^s_{\Box}:=\mathbf{V}(\text{div};\Box) 
\cap \mathbf{H}^s_0(\Box)$, $s \ge 0$, and its biorthogonal dual $\tilde {\bf \Psi}^\Box=\{{\tilde \psi}^\Box_{j,k}\}_{j 
\geq j_0, k \in \mathcal{J}_j^\Box}$. It holds that ${\bf \Psi}^\Box$ is a frame for 
$\mathbf{V}^s_{\Box}$.
Our aim is to show that the system
\begin{equation}
\label{aggframe}
{\bf \Psi} := ( \psi_{j,k}^i)_{(i,j,k)\in \Lambda},
\end{equation}
where
\begin{equation}
\label{lift}
\psi_{j,k}^i(x):=\frac{ \nabla \kappa_i\bigl(\kappa_i^{-1}(x)\bigr) \cdot \psi^{\Box}_{j,k}\bigl(\kappa_i^{-1}(x)\bigr)}{
\left|\det \nabla \kappa_i\bigl(\kappa_i^{-1}(x)\bigr)\right|^{1/2}},\quad
\mbox{for all } i=1,...,M, \ j \geq j_0, \ k \in \mathcal{J}^\Box_j, \ x\in\Omega_i,
\end{equation}
and $ \psi_{j,k}^i(x) = \mathbf{0}$ for $x\in\Omega\setminus\Omega_i$,
is a frame for $\mathbf{V}^s:= \mathbf{V}(\text{div};\Omega) \cap \mathbf{H}^s_0(\Omega)$. 
Anagously, we define its local duals by
\begin{equation}
\label{lift2}
\tilde \psi_{j,k}^i(x):=\frac{\nabla \kappa_i\bigl(\kappa_i^{-1}(x)\bigr) \cdot \tilde \psi^{\Box}_{j,k}\bigl(\kappa_i^{-1}(x)\bigr)}{
\left|\det \nabla \kappa_i\bigl(\kappa_i^{-1}(x)\bigr)\right|^{1/2}},\quad
\mbox{for all } i=1,...,M, \ j \geq j_0, \ k \in \mathcal{J}^\Box_j, \  x\in\Omega_i,
\end{equation}
and $\tilde \psi_{j,k}^i(x) = \mathbf{0}$ for $x\in\Omega\setminus\Omega_i$.

First of all observe that under assumption (C1) each ${\bf \Psi}^i:=( \psi_{j,k}^i)_{j\geq j_0, k \in \mathcal{J}^\Box_j}$ is again a divergence-free wavelet frame for $\mathbf{V}_i^s:=\mathbf{V}(\text{div};\Omega_i) \cap \mathbf{H}^s_0(\Omega_i)$. 
Moreover, by \cite{U2,U1,U}  for $\uu_i = \sum_{j,k} c_{j,k} \psi_{j,k}^i$ the 
following   norm equivalence holds
\begin{equation}
\label{equiv}
\| \uu_i\|_{\mathbf{H}^s} \sim \left( \sum_{j,k} 2^{2s j} | c_{j,k}|^2 \right)^{1/2}, \quad s \ge 0.
\end{equation}

We next  prove the following two auxiliary lemmas.

\begin{lemma}
\label{boundedcoef}
Let ${\bf \Psi}^{\Box,1}$ and ${\bf \Psi}^{\Box,2}$ be two generic wavelet bases on $\Box$ with sufficient 
regularity and number of vanishing moments. We do not necessarily require that these wavelets vanish on 
the boundary. Moreover, let
$\Omega$ be affine (or conformal if $n=2$) image of $\Box$. 
Let the corresponding wavelet systems ${\bf \Psi}^1=(\psi_{j,k}^{1})_{j\geq j_0,k \in \mathcal{J}^{\Box,1}_j}$ and 
${\bf \Psi}^1=(\psi_{j,k}^{2})_{j\geq j_0,k \in \mathcal{J}^{\Box,2}_j}$ on $\Omega$ be constructed by lifting as in 
\eqref{lift}. Then the cross Gramian matrix
\begin{equation}
\label{crossGram}
	M^{{\bf \Psi}^1,{\bf \Psi}^2}:= (\langle   \psi_{j,k}^{1}, \psi_{j',k'}^{2} \rangle )_{(j,k) \in \Lambda^1, (j',k') \in \Lambda^2}
\end{equation}
defines a bounded operator from $\ell_{2,2^{s j}}(\Lambda^2)$ onto $\ell_{2,2^{s j}}(\Lambda^1)$ with
$$
	\ell_{2,2^{s j}}(\Lambda) := \left \{ (c_{j,k})_{(j,k) \in \Lambda}: \left (\sum_{j \geq j_0, k 
	\in \mathcal{J}^{\Box}_j} 2^{2 s j} |c_{j,k}|^2 \right)^{1/2}< \infty \right \},
$$
for any $|s| \leq s_0$, where $s_0 \in \mathbb{R}_+$ depends on the regularity and the number of vanishing 
moments of ${\bf \Psi}^{\Box,1}$ and ${\bf \Psi}^{\Box,2}$.
\end{lemma}
\begin{proof}
The result is well--known and it can be immediately achieved, e.g., by combining 
\cite[Theorem 6.5]{DFR} and \cite[Proposition 8.2]{DFR}, where in the latter one replaces
$\rho_1$ with $\rho_2$.
\end{proof}

\begin{lemma}
\label{lastlem}
Let
$\cC:=\{\Omega_i\}_{i=1}^M$ be an overlapping, relatively compact covering of
$\Omega$ satisfying (C1).
Let the projections of $\LL_2(\Omega)$ onto $\bV_i:=\mathbf{V}(\emph{div},\Omega_i)$ be
\begin{equation}
P_{\mathbf{V}_i}(\uu):= \sum_{j \geq j_0,k \in \mathcal{J}^{\Box}_j} \langle \uu ,\tilde \psi_{j,k}^i 
\rangle \psi_{j,k}^i, \quad i=1,\dots, M.
\end{equation}
Then
\begin{itemize}
\item[(i)] $\displaystyle{\left ( \sum_{j \geq j_0,k \in \mathcal{J}^{\Box}_j} 2^{2sj} |\langle \uu, 
\tilde \psi_{j,k}^i \rangle |^2 \right )^{1/2} \lesssim  \|\uu \|_{\mathbf{H}^s}}$ for any
$i=1,\dots, M$ and all $\uu \in \mathbf{H}^s(\Omega)$;
\item[(ii)] $\| P_{\mathbf{V}_i}(\uu) \|_{\mathbf{H}^s} \lesssim \|\uu \|_{\mathbf{H}^s}$
 for any $i=1,\dots,M$ and all $\uu \in \mathbf{H}^s(\Omega)$.
\end{itemize}
\end{lemma}

\begin{proof}
Fix $i=1,\dots, M$. Cleary $ \uu|_{\Omega_i}  \in \mathbf{H}^s(\Omega_i)$, if $\uu \in \mathbf{H}^s(\Omega)$. 
It is well--known that $\uu|_{\Omega_i} = \sum_{j \geq j_0, k\in \mathcal{J}^{\Box,\circ}_j} \langle 
\uu, \tilde \psi_{j,k}^{i,\circ} \rangle \psi_{j,k}^{i,\circ}$ for some suitable wavelet basis 
$\{\psi_{j,k}^{i,\circ}\}_{j \geq j_0, k\in \mathcal{J}^{\Box,\circ}_j}$
constructed by lifting  a given reference basis on $\Box$. Moreover,  
\begin{equation}
\label{refnorm}
\left ( \sum_{j \geq j_0,k \in \mathcal{J}^{\Box,\circ}_j} 2^{2sj} |\langle \uu, 
\tilde \psi_{j,k}^{i,\circ} \rangle |^2 \right )^{1/2} \lesssim  \|\uu \|_{\mathbf{H}^s}	
\end{equation}
and
$$
\langle \uu, \tilde \psi_{j',k'}^i \rangle = 
\sum_{j \geq j_0, k\in \mathcal{J}^{\Box,\circ}_j} \langle \uu, 
\tilde \psi_{j,k}^{i,\circ} \rangle \langle \psi_{j,k}^{i,\circ}, \tilde \psi_{j',k'}^i \rangle,
\quad j' \ge j_0, \ k \in J_{j'}^{\Box}.
$$
Combining \eqref{refnorm} and Lemma \ref{boundedcoef} we get (i). Due to \eqref{equiv} we have
$$
\| P_{\mathbf{V}_i}(\uu) \|_{\mathbf{H}^s} \sim  \left( \sum_{j' \ge j_0, k' \in \mathcal{J}_{j'}^{\Box}} 2^{2s j'} | \langle \uu, \tilde 
\psi_{j',k'}^i \rangle|^2 \right)^{1/2} \lesssim \left ( \sum_{j \geq j_0,k \in \mathcal{J}^{\Box,
\circ}_j} 2^{2sj} |\langle \uu, \tilde \psi_{j,k}^{i,\circ} \rangle |^2 \right )^{1/2} \lesssim  
\|\uu \|_{\mathbf{H}^s}.
$$
\end{proof}

We next present the frame construction that allows us to work with $\Omega$, which are not
necessarily affine images of $(0,1)^n$, and provides the
characterization of $\bV^s$ for any $s \geq 1$.  We require that $\cC$ is such that
\begin{itemize}
\item [(C2)] it holds
\begin{equation} \label{C2}
P_{V_i} (\vv)|_{\Omega_i \setminus \cup_{k=2}^{\tilde{M}} \Omega_{i_k}} \equiv \vv|_{\Omega_i \setminus 
\cup_{k=2}^{\tilde{M}} \Omega_{i_k}},
\end{equation}
for any $\vv \in \bV^s$ such that $\vv|_{\partial \text{supp}(\vv)} = 0$,
$\Omega_i \subseteq \text{supp}(\vv) \subseteq \cup_{k=1}^{\tilde{M}} \Omega_{i_k}$ 
for $\{i_1, ..., i_{\tilde{M}}\} 
\subseteq \{1, ..., M\}$, without loss of generality $i:=i_1$. 
\end{itemize}

\begin{rems}

1. Note that the domain decomposition satisfying (C2) is necessarily overlapping. If it were
not the case, then it would hold that $\Omega_i=\Omega_i \setminus \cup_{\{1,\dots ,M\} \setminus \{i\}}
\Omega_j$ and, by \eqref{C2}, $\mathbf{P}_{V_i} \vv=\vv|_{\Omega_i}$. This, would imply that 
$\vv|_{\partial \Omega_i}=\mathbf{P}_{V_i} \vv|_{\partial \Omega_i}=0$, which is not necessarily the case
 for an arbitrary $\vv \in \bV^s$.
 
2. To illustrate our frame construction we consider $\Omega = \Omega_1 \cup \Omega_2$, which is not an affine image of $(0,1)^2$,
as in Figure 5. It holds $\mathbf{P}_{V_1} \vv=\vv$ on $\Omega \setminus \Omega_2$ and $\mathbf{P}_{V_1} \vv 
\in \bV^s$, $s \ge 1$. This, implies that $\vv-\mathbf{P}_{V_1} \vv \in \bV^s$. Note that by construction 
$\text{supp} (\vv-\mathbf{P}_{V_1} \vv) \subseteq \Omega_2$ and $\vv-\mathbf{P}_{V_1} \vv|_{\partial \Omega_2}=0$. 
Therefore, we obtain the following decomposition 
$\vv=\mathbf{P}_{V_1} \vv+\mathbf{P}_{V_2}(\vv-\mathbf{P}_{V_1} \vv)$.

\begin{figure}[ht]
\hbox to \hsize {\hfill \epsfig{file=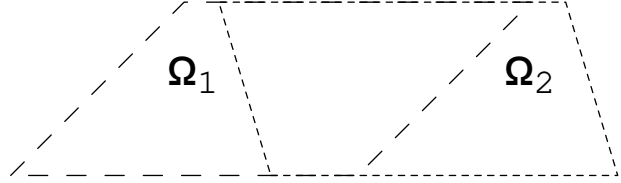,width=8.5cm} \hfill}
\caption{An example of Overlapping Domain Decomposition satisfying (C1) and (C2). }
\end{figure}
\end{rems}

The following Theorem generalizes the argument in Remark 2. above to domain decompositions 
consisting of more than two patches.

\begin{tm}
\label{decomp2}
Let $\cC:=\{\Omega_i\}_{i=1}^M$ be an overlapping, relatively compact covering of
$\Omega$ satisfying (C1)--(C2). Let   $\Psi := ( \psi_{j,k}^i)_{(i,j,k)\in \Lambda}$ be as in \eqref{aggframe}
and $P_{\mathbf{V}_i}$ as in Lemma \ref{lastlem}.
For $\uu \in \bV^s$, $s \ge 1$,  consider the functions and domains defined by the following recursion 
$$
\uu^{(0)} = \uu, \quad \Omega^{(0)} = \Omega,
$$
\begin{equation}
\label{induc}
\left\{
\begin{array}{l} 
\uu_{n+1} = P_{V_{n+1}} \uu^{(n)}\\
\uu^{(n+1)} = \uu^{(n)} - \uu_{n+1}\\
\Omega^{(n+1)} = \emph{supp}(\uu^{(n)})
\end{array}\right . , \  n=0, \dots, M-1.
\end{equation}
Then
\begin{eqnarray}
\label{framedecfin}
\uu = \sum_{i=1}^M \uu_{i} 
=\sum_{i=1}^M \sum_{j\geq j_0, k \in \mathcal{J}^\Box_j} \langle \uu^{(i-1)}, \tilde \psi_{j,k}^i \rangle \psi_{j,k}^i
\end{eqnarray}
with
\begin{equation}
\label{normeq2}
\|\uu\|_{\mathbf{H}^s} \sim  \left( \sum_{i=1}^M  \sum_{j\geq j_0, k \in \mathcal{J}^\Box_j} 2^{2 s j}|\langle \uu^{(i-1)}, \tilde \psi_{j,k}^i \rangle |^2 \right)^{1/2}.
\end{equation}
In particular, ${\bf \Psi}$ is an $\LL_2$ frame for $\bV_0(\emph{div}; \Omega)$ and a Gelfand frame for $(\bV^s, \bV_0(\emph{div}; \Omega),(\bV^s)')$.
\end{tm}

\begin{proof}
By property (C2), Lemma \ref{lastlem} and by induction we get that $\uu^{(n)} \in \bV^s$ with $ \|\uu^{(n)} 
\|_{\bH^s} \lesssim \|\uu \|_{\bH^s}$, $\uu^{(n)}|_{\partial \Omega^{(n+1)}} = 0$ and $\Omega^{(n)} \subseteq \cup_{i=n+1}^M 
\Omega_i$. Thus, $\Omega^{(M)} = \emptyset$ and, therefore, in the recursive definition \eqref{induc} 
we have $0 \leq n \leq M-1$. This implies also that  $\text{supp}(\uu^{(M-1)}) \subseteq \Omega_M$ and $\uu_{M} =  
P_{V_M} \uu^{(M-1)} = \uu^{(M-1)}$. And by induction we gets  
$\uu^{(M-1)} = \uu- \sum_{i=1}^{M-1} \uu_i$. This immediately implies 
$$
\uu =  \sum_{i=1}^{M-1} \uu_{i}+ \left (\uu - \sum_{i=1}^{M-1} \uu_{i} \right ) =\sum_{i=1}^{M-1} \uu_{i}+ \uu^{(M-1)} = \sum_{i=1}^M \uu_{i}
$$
and, consequently, the wavelet decomposition of $\uu$
$$
\uu =  \sum_{i=1}^M \uu_{i} = \sum_{i=1}^M P_{\bV_i} \uu^{(i-1)} = \sum_{i=1}^M \sum_{j\geq j_0, k \in \mathcal{J}^\Box_j} \langle \uu^{(i-1)}, \tilde \psi_{j,k}^i \rangle \psi_{j,k}^i.
$$
By an application of Lemma \ref{lastlem} (i) and property \eqref{equiv} one has
$$
\|\uu\|_{\mathbf{H}^s} \lesssim  \left( \sum_{i=1}^M \sum_{j\geq j_0, k \in \mathcal{J}^\Box_j}  2^{2 s j} |\langle \uu^{(i-1)}, \tilde \psi_{j,k}^i \rangle |^2 \right)^{1/2}.
$$
Conversely, by Lemma \ref{lastlem} (i)
$$
\left( \sum_{i=1}^M \sum_{j\geq j_0, k \in \mathcal{J}^\Box_j}  2^{2 s j} |\langle \uu^{(i-1)}, 
\tilde \psi_{j,k}^i \rangle |^2 \right)^{1/2} \lesssim \left( \sum_{i=1}^M \|\uu^{(i-1)}\|^2_{\bH^s} \right )^{1/2} \lesssim \|\uu \|_{\bH^s}.
$$
Note that the norm equivalences above hold for also for $s=0$ and, thus, 
by density of $\bV^s$ in $\bV_0(\text{div}; \Omega)$, ${\bf \Psi}$ is an $\LL_2$ 
frame for $\bV_0(\text{div}; \Omega)$ and a Gelfand frame for $(\bV^s, \bV_0(\text{div}; 
\Omega),(\bV^s)')$.
\end{proof}

\begin{rems}
1. As mentioned in Remark 2. after Theorem \ref{cmpxtm}, the wavelet frame expansion \eqref{framedecfin} allows for the construction of an admissible projector $\tilde \PP$ for the optimal application of the procedure $\mathbf{modSOLVE}$. Such construction follows the one in 
\cite[Section 4.4]{S}.

2. The aggregate wavelet frames constructed  in  
Theorem \ref{decomp2}  are multiscale and  produce 
coefficients that can be naturally structured into suitable trees, as the ones in \cite{DS1,DS2}. Thus, 
the algorithms defined in \cite{CDD4}, see (A3), for 
the evaluation of nonlinear functionals can also be applied in our case. 
The optimality  also requires the characterization of Besov spaces in terms of norm equivalences. 
As soon as the local divergence--free bases on $\Box$ provides such a characterization, one can show similarly 
to Theorem \ref{decomp2}  that the global divergence--free basis provides the  characterization of Besov spaces
on the  whole domain. We postpone the details to a forthcoming paper. 

3. In practice, one  never  uses the full frame ${\bf \Psi}:=\cup_{i=1}^M {\bf \Psi}_i = 
\{ {\psi}_{j,k}^i\}_{j \geq -1, k \in 
\mathcal{J}_j^i, i=1,...,M}$, but rather ${\bf \Psi}_J:=\{ {\psi}_{j,k}^i\}_{-1 \leq j 
\leq J, k \in \mathcal{J}_j^i, i=1,...,M}$ (for example, in the concrete implementation of the schemes discussed 
in \cite{W}, an upper bound in terms of concretely used scale levels $j$ must be enforced for computational 
(memory/storage) limits, without spoiling the accuracy). Therefore, the first approach mentioned above is 
``ready--to--use''. 
It is well--known that any finite system is already  a (Gelfand) frame for its span.  Using this fact one 
can consider  solving {\bf Problem 4}  
 in $\text{span}({\bf \Psi}_J) \subseteq \mathbf{V}(\text{div};\Omega)\cap \mathbf{H}^1_0(\Omega)$ for $J\geq 0$ large. 
 The scheme we propose is  well--defined in this case.
\end{rems}

\section{Conclusion}

We  present the  first {\it convergent and implementable adaptive scheme} based on frame decompositions 
for the numerical integration of magnetohydrodynamic flows. 
Certainly, it is impossible to construct a universal scheme suitable for any physical problem.
Thus, the convergence of the algorithm is ensured under certain assumptions on the physical parameters. Such
 assumptions are standard when dealing with nonlinear problems of fluid dynamics and magnetohydrodynamics 
and still lead to the analysis of physically realistic problems, see \cite{MS1,MS2,MS3}. 

We  show that under Dirichlet boundary conditions on the velocity and 
the zero flux boundary conditions on the electric current the velocity and electric current are divergence--free vector valued functions. 
The discretization of the problem is, therefore, realized by expanding the solution with respect to certain 
divergence--free wavelet frames constructed on Overlapping Domain Decompositions. 
The construction  based on ODD is definitively 
a breakthrough in comparison with \cite{U2,U1,U}. There divergence--free wavelet bases 
have been defined in a relatively simple way essentially only for domains that are affine images of cubes. 
Our construction of frames allows for rather general polygonal domains and avoids the use of fictitious domain
techniques, see, e.g., \cite{DKlU}.  
These frames can be also used for the discretization of Stokes and Navier--Stokes equations \cite{DVU}.
Of importance is that these aggregate divergence--free wavelet frames preserve 
multiscale properties of standard wavelet bases and their capability to characterize certain function spaces. 
These properties ensure the optimal convergence 
rates and complexity of the scheme we propose, if the solutions are in certain sparseness 
classes of functions. It has not yet been theoretically proven that the solutions belong to such function spaces (e.g., 
Besov spaces), numerical simulations \cite{U3} though  support this assumption.

\section*{Acknowledgement} 
M. Fornasier acknowledges the financial support provided through the
Intra-European Individual Marie Curie Fellowship Programme, 
under contract MEIF-CT-2004-501018. All of the authors acknowledge 
the hospitality of Dipartimento di Metodi e Modelli Matematici 
per le Scienze Applicate, Universit\`a di Roma ``La Sapienza'', 
Italy, during the preparation of this work.
The authors want to thank Daniele Boffi and Karsten Urban for the helpful and fruitful discussions on divergence-free function spaces.

\bibliography{banach}

\newpage


\noindent Maria Charina\\
Universit\"at Dortmund\\
Fachbereich Mathematik, Lehrstuhl VIII\\
Vogelpothsweg 87\\
44221 Dortmund\\
Germany\\
email: {\tt maria.charina@math.uni-dortmund.de}\\
WWW: {\tt http://www.mathematik.uni-dortmund.de/lsviii/charina/Engindex.html}
\\

\noindent Costanza Conti\\
Universit\`a di Firenze\\
Dipartimento di Energetica ``Sergio Stecco''\\
Via C. Lombroso 6/17\\
50134 Firenze\\
Italy\\
email: {\tt c.conti@ing.unifi.it}\\
WWW: {\tt http://www.de.unifi.it/anum/Conti/conti$\_$En.htm}
\\

\noindent Massimo Fornasier\\
Universit\`a ``La Sapienza'' in Roma\\
Dipartimento di Metodi e Modelli Matematici per le Scienze Applicate\\
Via A. Scarpa, 16/B\\
00161 Roma\\
Italy\\
and\\
Universit\"at Wien\\
NuHAG, Fakult\"at Mathematik\\
Nordbergstrasse 15\\
1090 Wien\\
Austria\\
email: {\tt mfornasi@math.unipd.it}\\
WWW: {\tt http://www.math.unipd.it/$\sim$mfornasi}

\end{document}